\documentstyle[11pt,amssymb]{article}

\begin{document}

\begin{center}
{\large {\bf Riemann-Roch isomorphism, Chern-Simons invariant}} 
\\
{\large {\bf and Liouville action}} 
\end{center}

\begin{center}
By {\sc Takashi Ichikawa} 
\end{center}

\begin{center}
{\bf Abstract} 
\end{center}

Using the arithmetic Schottky uniformization theory, 
we show the arithmeticity of $PSL_{2}({\mathbb C})$ Chern-Simons invariant. 
In terms of this invariant, we give an explicit formula of the Riemann-Roch isomorphism 
as Zograf-Mcintyre-Takhtajan's infinite product for families of algebraic curves. 
By this formula, we determine the unknown constant 
which appears in the holomorphic factorization formula of  
determinant of Laplacians on Riemann surfaces via the classical Liouville action. 
As an application, we show the rationality of Ruelle zeta values for 
Schottky uniformized $3$-manifolds. 
\vspace{2ex}

\begin{center}
{\bf 1. Introduction} 
\end{center}

Arakelov theory, Chern-Simons theory and Liouville theory 
have different origins and important roles in many areas of mathematics 
containing arithmetic (algebraic) geometry, differential geometry, 
global analysis and mathematical physics. 
The aim of this paper is to show that a combination of the $3$ theories 
contributes to these advances. 
More precisely, 
we consider together the subjects: 
\begin{itemize}

\item 
the arithmetic Riemann-Roch theorem for families of algebraic curves 
(cf. [D2, Fr, GS, W]), 

\item 
the $PSL_{2}({\mathbb C})$ Chern-Simons theory for 
hyperbolic $3$-manifolds with boundary (cf. [GM, MP]), 

\item 
the holomorphic factorization formula of determinant of Laplacians 
on Riemann surfaces in terms of the classical Liouville action 
(cf. [Z1, Z2, MT]). 

\end{itemize}
As a result of this consideration, 
we have the following results: 
\begin{itemize}

\item 
arithmeticity of the power series expansion of 
the $PSL_{2}({\mathbb C})$ Chern-Simons invariant (cf. Theorem 5.1), 

\item 
an explicit formula of the Riemann-Roch isomorphism 
under the trivialization by the arithmetic $PSL_{2}({\mathbb C})$ 
Chern-Simons invariant and Eichler cohomology (cf. Theorem 5.5), 

\item 
determination of the unknown constant in the holomorphic factorization formula of 
Zograf and Mcintyre-Takhtajan (cf. Theorems 5.7 and 5.8). 

\end{itemize}

We will review the main idea of this paper. 
The arithmetic Riemann-Roch theorem is an advanced, 
i.e., metrized version of the Grothendieck Riemann-Roch theorem, 
and is especially applied to Diophantine problem. 
For a family of algebraic curves, 
these theorems gives an isometric isomorphism (up to a nonzero constant) 
$$
\lambda_{k}^{\otimes 12} \cong \kappa^{\otimes d_{k}}; \ 
d_{k} := 6 k^{2} - 6 k + 1, 
$$ 
where $\lambda_{k}$ denotes the $k$th tautological line bundle 
with Quillen metric, 
and $\kappa$ denotes the Deligne pairing of the relative canonical sheaf 
with itself. 
The Liouville (field) theory gives rise to the Liouville action 
which is a functional on the space of conformal metrics on Riemann surfaces 
whose critical points correspond to hyperbolic metrics. 
Then the holomorphic factorization formulae relate the classical Liouville action 
with the determinants of Laplacians 
in terms of Zograf-Mcintyre-Takhtajan's infinite products 
which are extensions of Ramanujan's delta function. 
Notice that these determinants corresponds to 
the Quillen metric on $\lambda_{k}$, 
and as is shown by Aldrovandi [A], 
the classical Liouville action gives the natural metric on $\kappa$. 
Therefore, 
one can suppose that the arithmetic Riemann-Roch isometry is equivalent to 
the holomorphic factorization formula. 

In order to make this equivalence more complete, 
we use the theories of arithmetic Schottky uniformization [I2] 
and of Chern-Simons line bundles [GM]. 
In [I2], one has a higher genus version of the Tate curve, 
called {\it generalized Tate curves}, 
which become Schottky uniformized Riemann surfaces over ${\mathbb C}$  
and give local coordinates on the moduli space over ${\mathbb Z}$ of algebraic curves. 
Therefore, based on the ${\mathbb Z}$-structure of $\lambda_{k}$ and $\kappa$, 
we will describe the Riemann-Roch isomorphism using the following simple fact: 
{\it if two power series over ${\mathbb Z}$ are known to be proportional 
to each other and primitive (i.e., not congruent to $0$ modulo any prime), 
then these power series are equal up to a sign.} 

In [GM], it is shown that the $PSL_{2}({\mathbb C})$ Chern-Simons invariants 
of Schottky uniformized $3$-manifolds give rise to a hermitian line bundle, 
called the {\it Chern-Simons line bundle}, 
${\cal L}$ which is isometrically isomorphic to $\lambda_{1}^{\otimes (-6)}$ 
on the Schottky space. 
We will show that ${\cal L}^{\otimes 2}$ is actually defined 
on the moduli space of Riemann surfaces, 
and that there exists an isometric isomorphism (up to a nonzero constant) 
$$
\lambda_{1}^{\otimes 12} \cong {\cal L}^{\otimes (-2)}. 
$$
By the above isomorphisms, 
we have an identification 
$$
{\cal L}^{\otimes 2} \cong \kappa^{\otimes (-1)}
$$ 
as line bundles over the Deligne-Mumford compactification of the moduli space. 
Then by the arithmetic Schottky uniformization theory, 
one can modify the exponential of the $PSL_{2}({\mathbb C})$ Chern-Simons invariant 
such as to have a universal expression as a power series over ${\mathbb Z}$ 
which we call the {\it arithmetic Chern-Simons invariant}. 
Therefore, this gives a (local) canonical trivialization of  $\kappa$. 
Furthermore, we show that the natural ${\mathbb Z}$-structure of 
the Eichler cohomology group gives a canonical trivialization of $\lambda_{k}$  
by studying the behavior under the degeneration of a generalized Tate curve 
to a pointed Tate curve. 
Therefore, by the holomorphic factorization formula and the above simple fact, 
we have: 
\vspace{2ex}

{\sc Theorem 1.1} (see Theorem 5.5 for the precise statement). 
\begin{it} 
Let $F_{k}$ denote Zograf-Mcintyre-Takhtajan's infinite product whose main part is 
$$
\prod_{\{ \gamma \}} \prod_{m=k}^{\infty} \left( 1 - q_{\gamma}^{m} \right), 
$$
where $\{ \gamma \}$ runs over primitive conjugacy classes in a Schottky group $\Gamma$,  
and $q_{\gamma}$ denotes the multiplier of $\gamma$. 
Then under the trivializations of $\kappa$ and $\lambda_{k}$ 
by the arithmetic Chern-Simons invariant and the Eichler cohomology group  respectively, 
the Riemann-Roch isomorphism 
$\lambda_{k}^{\otimes 12} \stackrel{\sim}{\rightarrow} \kappa^{\otimes d_{k}}$ 
for families of algebraic curves is expressed as $\pm F_{k}^{-12}$. 
\end{it} 
\vspace{2ex}

Using this theorem and the arithmetic Riemann-Roch theorem, 
we also determine the unknown constant in the holomorphic factorization formula. 

We give some comments on related works. 
In [MP, Theorem 1.1], 
Mcintyre and Park express the $PSL_{2}({\mathbb C})$ Chern-Simons invariants 
of Schottky uniformized $3$-manifolds 
in terms of the Bergman tau function (cf. [KK1]) and Zograf's function. 
Since this tau function has the modular property (cf. [KK2, (3.15)]), 
Mcintyre-Park's formula seems to represent the isomorphism 
$\lambda_{1}^{\otimes 24} \cong {\cal L}^{\otimes (-4)}$. 
Another relationship between the $PSL_{2}({\mathbb C})$ Chern-Simons invariants 
and Deligne pairings are given by Morishita and Terashima [MoT] 
for $3$-manifolds obtained as knot complements. 

As an application of the above results, 
we consider the rationality of special values of the Ruelle zeta function 
$$
Z_{\Gamma}(s) = \prod_{\{ \gamma \}} \left( 1 - |q_{\gamma}|^{s} \right)^{-1}, 
$$
for a Schottky group $\Gamma$ 
based on the fact that if $\Gamma \subset PSL_{2}({\mathbb R})$, 
then Zograf-McIntyre-Takhtajan's infinite products become Selberg type zeta values. 
A conjecture of Deligne [D1], 
which is extended by Beilinson [B] and Bloch-Kato [BK], 
states the following: 
each special value of the Hasse-Weil $L$-function of a motive ${\cal M}$ 
becomes the product of its transcendental part represented by 
the period or regulator of ${\cal M}$, 
and its rational part represented by the arithmetic invariant of ${\cal M}$. 
As its analogy for {\it geometric} zeta functions, 
we show: 
\vspace{2ex}

{\sc Theorem 1.2} (see Theorem 6.2 and Corollary 6.3 for the precise statement). 
\begin{it} 
Assume that $\Gamma \subset PSL_{2}({\mathbb R})$, 
and that the Riemann surface uniformized by $\Gamma$ has a model $C_{\Gamma}$ 
as an algebraic curve defined over a subfield $K$ of ${\mathbb C}$. 
Then the modified Ruelle zeta value of $\Gamma$ at any integer $k > 1$ 
belongs to 
$$
\frac{\mbox{period of regular $(k+1)$-forms on $C_{\Gamma}$}}
{\mbox{period of regular $k$-forms on $C_{\Gamma}$}} 
\cdot c(\Gamma)^{12k} \cdot K^{\times}, 
$$
where $c(\Gamma)$ denotes the ``value'' at $\Gamma$ of Zograf's infinite product 
divided by the period of regular $1$-forms on $C_{\Gamma}$. 
\end{it} 
\vspace{2ex}

We also express this zeta value by the discriminant of $C_{\Gamma}$ 
based on results of Saito [S], 
however we have no result on the rational or transcendental property of 
$c(\Gamma)$. 
Notice that using a result of [I3], 
one can construct Schottky groups $\Gamma \subset PSL_{2}({\mathbb R})$ 
such that $C_{\Gamma}$ are hyperelliptic and defined over ${\mathbb Q}$. 

By using the Selberg trace formula, 
the leading term at the central point $0$ of the Ruelle $L$ function 
of a hyperbolic manifold $X$ was studied by Fried [Fri] when $X$ is compact, 
and by Park [P], Sugiyama [Su1-5] and Gon-Park [GP] when $X$ has finite volume. 
Especially, Sugiyama's results are regarded as a geometric analog to Iwasawa theory 
and the Beilinson conjecture. 
In our case where $C_{\Gamma} \otimes_{K} {\mathbb C}$ 
is the boundary of the hyperbolic $3$-manifold $X$ uniformized by $\Gamma$, 
the above theorem states that the rationality of the Ruelle zeta values of $X$ 
at positive integers is controlled by the field of definition of its boundary. 
\vspace{2ex}

{\it Acknowledgments}. 
The author would like to thank deeply J. Park for his introduction 
of Zograf's formula which gave a motivation of this work. 
The author also thanks A. McIntyre, K. Sugiyama and S. Koyama 
for their valuable comments. 
This work is partially supported by the JSPS Grant-in-Aid for
Scientific Research (C) 26400018. 
\vspace{2ex}

\begin{center}
2. {\bf  The Chern-Simons line bundle} 
\end{center}

2.1. {\it $PSL_{2}({\mathbb C})$ Chern-Simons invariant}. 
For a Riemannian oriented $3$-manifold $X$ with metric $g$ and 
an orthonormal frame $S$ over $X$, 
Chern and Simons [CS] introduced the Chern-Simons invariant ${\rm CS}(g, S)$ 
which is defined as  
$$
\frac{1}{16 \pi^{2}} \int_{X} S^{*} \left( {\rm Tr} \left( 
\omega \wedge d \omega 
+ \frac{2}{3} \omega \wedge \omega \wedge \omega \right) \right),  
$$
where $\omega$ denotes the Levi-Civita connection form for $g$. 
This complexified version, called the $PSL_{2}({\mathbb C})$ Chern-Simons invariant, 
was given by Yoshida [Y]. 
Furthermore, Guillarmou and Moroianu [GM], 
McIntyre and Park [MP] extended this invariant for nonclosed  $3$-manifolds, 
and they studied the relationship between the associated hermitian holomorphic 
line bundle and Quillen's determinant line bundle in the Schottky setting. 
Recall that in [GM, Proposition 16] and [MP, 4.5], 
for each Schottky uniformized $3$-manifold $X$ with hyperbolic metric $g$, 
the $PSL_{2}({\mathbb C})$ Chern-Simons invariant 
${\rm CS}^{PSL_{2}({\mathbb C})}$ is defined and expressed as 
$$
{\rm CS}^{PSL_{2}({\mathbb C})} (g, S) 
= - \frac{\sqrt{-1}}{2 \pi^{2}} {\rm Vol}_{\rm R}(X) 
+ \frac{\sqrt{-1}}{2 \pi} \chi (\partial X) + {\rm CS}(g, S). 
$$
Here ${\rm Vol}_{\rm R}(X)$ denotes the renormalized volume of $X$, 
and $\chi (\partial X)$ denotes the Euler characteristic of the boundary of $X$.  
\vspace{2ex}

2.2. {\it Schottky uniformization}. 
Schottky groups of rank $g$ are free groups generated by 
$\gamma_{1},..., \gamma_{g} \in PSL_{2}({\mathbb C})$ which map Jordan curves  
$C_{1},... , C_{g} \subset 
{\mathbb P}^{1}({\mathbb C}) = {\mathbb C} \cup \{ \infty \}$ 
to other Jordan curves 
$C_{-1},... , C_{-g} \subset {\mathbb P}^{1}({\mathbb C})$ 
(with orientation reversed). 
Each element $\gamma \in \Gamma - \{ 1 \}$ is conjugated in 
$PSL_{2}({\mathbb C})$ to $z \mapsto q_{\gamma} z$ 
for some $q_{\gamma} \in {\mathbb C}^{\times}$ with $|q_{\gamma}| < 1$, 
called the {\it multiplier} of $\gamma$. 
Therefore, 
$$
\frac{\gamma(z) - a_{\gamma}}{\gamma(z) - b_{\gamma}} = 
q_{\gamma} \frac{z - a_{\gamma}}{z - b_{\gamma}} 
$$
for some element $a_{\gamma}, b_{\gamma}$ of ${\mathbb P}^{1}({\mathbb C})$ 
called the {\it attractive}, {\it repulsive} fixed points of $\gamma$ 
respectively. 
Then the discontinuity set 
$\Omega_{\Gamma} \subset {\mathbb P}^{1}({\mathbb C})$ 
under the action of $\Gamma$ has a fundamental domain $D_{\Gamma}$ 
which is given by the complement of the union of the interiors 
of $C_{i}$ $(i = \pm 1,..., \pm g)$. 
The quotient space $\Omega_{\Gamma}/\Gamma$ 
of $\Omega_{\Gamma}$ by $\Gamma$ is a (compact) Riemann surface 
of genus $g$ which we denote by $R_{\Gamma}$. 
Furthermore, by a result of Koebe, 
every Riemann surface of genus $g$ can be represented in this manner. 
A Schottky group $\Gamma$ is {\it marked} if its free generators 
$\gamma_{1},..., \gamma_{g}$ are fixed, 
and a marked Schottky group $(\Gamma; \gamma_{1},..., \gamma_{g})$ 
is {\it normalized} if 
$a_{\gamma_{1}} = 0$, $b_{\gamma_{1}} = \infty$ and $a_{\gamma_{2}} = 1$. 
By definition, the Schottky space ${\mathfrak S}_{g}$ of degree $g$ 
is the space of marked Schottky groups of rank $g$ 
modulo conjugation in $PSL_{2}({\mathbb C})$ which becomes 
the space of normalized Schottky groups of rank $g$ if $g > 1$. 
Then ${\mathfrak S}_{g}$ is a covering space of the moduli space of 
Riemann surfaces of genus $g$.
\vspace{2ex}

2.3. {\it Chern-Simons line bundle}. 
We define the Chern-Simons line bundle ${\cal L}_{{\mathfrak S}_{g}}$ 
over the Schottky space ${\mathfrak S}_{g}$ of degree $g$ 
following Freed [F], Ramadas-Singer-Weitsman [RSW] and especially 
Guillarmou and Moroianu [GM]. 

For each Schottky group $\Gamma \subset PSL_{2}({\mathbb C})$, 
denote by $X = X_{\Gamma} = {\mathbb H}^{3}/\Gamma$ 
the associated hyperbolic $3$-manifold, 
where $ {\mathbb H}^{3}$ denotes the $3$-dimensional hyperbolic space. 
Then the boundary of $X$ becomes a Riemann surface $R = R_{\Gamma}$ 
which is Schottky uniformized by $\Gamma$. 
Denote by $F(X)$ the frame bundle of $X$, 
and by $C^{\infty}_{\rm ext}(R, *)$ the space of sections in $C^{\infty}(R, *)$ 
which can be extended to $\overline{X} = X \cup R$. 
Let 
$$
c_{\Gamma} : C^{\infty}(R, F(X)) \times C^{\infty}_{\rm ext}(R, SO(3)) 
\rightarrow {\mathbb C}
$$ 
be the map defined as 
$$
c_{\Gamma} \left( \widehat{S}, a \right) := 
\exp \left( \frac{2 \pi \sqrt{-1}}{16 \pi^{2}} 
\left( \int_{R} {\rm Tr} (\widehat{\omega} \wedge da \wedge a^{-1}) + 
\int_{X} \frac{1}{3} {\rm Tr} 
\left( (\widetilde{a}^{-1} d \widetilde{a})^{3} \right) \right) \right), 
$$ 
where $\widehat{\omega}$ is the connection form of the Levi-Civita connection 
of the associated metric on $\overline{X}$ along $R$ in the frame $\widehat{S}$, 
and $\widetilde{a}$ is any smooth extension of $a$ on $\overline{X}$. 
Then $c_{\Gamma} \left( \widehat{S}, a \right)$ is well-defined by [GM, Lemma 20], 
and is seen to satisfy the cocycle condition 
$$
c_{\Gamma} \left( \widehat{S}, ab \right) = 
c_{\Gamma} \left( \widehat{S}, a \right) \cdot 
c_{\Gamma} \left( \widehat{S}a, b \right). 
$$
We define the complex vector space $L_{\Gamma}$ as the space of 
complex-valued functions $f$ on $C^{\infty}_{\rm ext}(R, F(X))$ which satisfy  
$$
f \left( \widehat{S}a \right) = 
c_{\Gamma} \left( \widehat{S}, a \right) \cdot f \left( \widehat{S} \right) \ 
\left( a \in C^{\infty}_{\rm ext}(R, SO(3)) \right). 
$$
By this condition, any element in $L_{\Gamma}$ is determined by its value 
on any frame which can be extended to $\overline{X}$, 
and hence one can define the Chern-Simons line bundle over ${\mathfrak S}_{g}$ as 
$$
{\cal L}_{{\mathfrak S}_{g}} := 
\bigsqcup_{\Gamma \in {\mathfrak S}_{g}} L_{\Gamma}. 
$$
Furthermore, 
the cocycle has absolute value $1$, 
and hence there exists a canonical hermitian metric 
$\langle \cdot, \cdot \rangle_{\rm CS}$ on ${\cal L}_{{\mathfrak S}_{g}}$ 
defined as 
$$
\langle f_{1}, f_{2} \rangle_{\rm CS} := 
f_{1} \left( \widehat{S} \right) \overline{f_{2} \left( \widehat{S} \right)} 
$$
for two sections $f_{1}, f_{2}$ of ${\cal L}_{{\mathfrak S}_{g}}$ and 
$\widehat{S} \in C^{\infty}(R, F(X))$. 
\vspace{2ex}

2.4. {\it Isomorphism with determinant line bundle}. 
We review a result of [GM] on an explicit isomorphism between 
the Chern-Simons line bundle with connection and hermitian structure 
and the determinant line bundle on the Schottky space. 

Let $R$ be a Riemann surface of genus $g > 0$, 
and $\left\{ \alpha_{1},..., \alpha_{g}, \beta_{1},..., \beta_{g} \right\}$ 
be a set of standard generators of $\pi_{1}(R, x_{0})$ for some $x_{0} \in R$ 
satisfying 
$$
\left( \alpha_{1} \beta_{1} \alpha_{1}^{-1} \beta_{1}^{-1} \right) \cdots 
\left( \alpha_{g} \beta_{g} \alpha_{g}^{-1} \beta_{g}^{-1} \right) = 1. 
$$ 
Then one can take a marked Schottky group 
$(\Gamma; \gamma_{1},..., \gamma_{g})$ such that $R = R_{\Gamma}$ 
and that each $C_{k}$ is homotopic to $\alpha_{k}$. 
Therefore, 
there is uniquely a basis $\varphi_{1},..., \varphi_{g}$ of holomorphic $1$-forms 
such that $\int_{\alpha_{j}} \varphi_{i}$ is equal to Kronecker's delta $\delta_{ij}$, 
and then the period matrix 
$\tau = \left( \int_{\beta_{j}} \varphi_{i} \right)$ becomes a symmetric matrix 
whose imaginary part is positive definite. 
For each Schottky group $\Gamma$ of rank $g$, 
$\{ \varphi_{1},..., \varphi_{g} \}$ is a basis of the space 
$H^{0} \left( R_{\Gamma}, \Omega_{R_{\Gamma}} \right)$ of 
holomorphic $1$-forms on $R_{\Gamma}$. 
Therefore, the {\it Hodge line bundle} $\lambda_{1}$ consisting of 
$\bigwedge^{g} H^{0} \left( R_{\Gamma}, \Omega_{R_{\Gamma}} \right)$ 
$(\Gamma \in {\mathfrak S}_{g})$ becomes a holomorphic line bundle 
on ${\mathfrak S}_{g}$ with a holomorphic canonical section 
$$
\varphi := \varphi_{1} \wedge \cdots \wedge \varphi_{g}. 
$$   
For each Riemann surface $R_{\Gamma}$, 
let $h$ be the associated hyperbolic metric, 
and $\det' \Delta_{h}$ be the (regularized) determinant of its Laplacian 
defined by Ray-Singer [RS]. 
Then the hermitian {\it Quillen metric} on $\lambda_{1}$ (cf. [Q]) is defined as 
$$
\left\| \varphi \right\|_{\rm Q}^{2} := 
\frac{\left\| \varphi \right\|_{h}^{2}}{\det' \Delta_{h}} = 
\frac{\det {\rm Im} (\tau)}{\det' \Delta_{h}} 
$$
at $\Gamma \in {\mathfrak S}_{g}$, 
where $\left\| \varphi \right\|_{h}$ is the hermitian product on 
$\bigwedge^{g} H^{0} \left( R_{\Gamma}, \Omega_{R_{\Gamma}} \right)$ induced by $h$. 
Therefore, there is the unique hermitian connection $\nabla^{\det}$ associated to 
the holomorphic structure on $\lambda_{1}$ and the hermitian norm 
$\left\| \varphi \right\|_{\rm Q}$. 

To describe the relationship between the Chern-Simons line bundle and 
the determinant line bundle, 
we use a formula proved by Zograf [Z1, Z2] whose generalization 
by McIntyre-Takhtajan [MT] will be reviewed below. 
Denote by $S_{\rm L} : {\mathfrak S}_{g} \rightarrow {\mathbb R}$ 
the classical Liouville action which is explicitly described by 
Takhtajan and Zograf [ZT] as 
\begin{eqnarray*} 
S_{\rm L} & = & 
\frac{\sqrt{-1}}{2} \int \int_{D_{\Gamma}} 
\left( \left| \frac{\partial \log \rho}{\partial z} \right|^{2} 
+ \rho \right) dz \wedge d \overline{z} 
\\ 
& & 
+ \ \frac{\sqrt{-1}}{2} \sum_{k=2}^{g} \oint_{C_{k}} 
\left( \log \rho - \frac{1}{2} \log \left| \gamma'_{k} \right|^{2} \right) 
\left( \frac{\gamma''_{k}}{\gamma'_{k}} dz - 
\frac{\overline{\gamma''_{k}}}{\overline{\gamma'_{k}}} d \overline{z} \right) 
\\
& & 
+ \ 4 \pi \sum_{k=2}^{g} \log \left| c(\gamma_{k}) \right|^{2}, 
\end{eqnarray*} 
where $D_{\Gamma}, C_{k} \subset {\mathbb P}^{1}({\mathbb C})$ are given in 2.2, 
$\rho(z) |dz|^{2}$ denotes the pullback of the hyperbolic metric on $R_{\Gamma}$ 
and $c(\gamma) = c$ for 
$\gamma = \left( \begin{array}{cc} a & b \\ c & d \end{array} \right)$. 
\vspace{2ex}

{\sc Theorem 2.1} (Zograf [Z1, Z2]). 
\begin{it} 
\begin{itemize}

\item[{\rm (1)}] 
There exists a holomorphic function 
$F_{1} : {\mathfrak S}_{g} \rightarrow {\mathbb C}$ 
such that
$$
\frac{\det' \Delta_{h}}{\det {\rm Im}(\tau)} = 
c_{g} \exp \left( - \frac{S_{\rm L}}{12 \pi} \right) |F_{1}(\Gamma)|^{2}, 
$$
where $c_{g}$ is a nonzero constant depending only on $g$. 

\item[{\rm (2)}] 
If the Hausdorff dimension $\delta_{\Gamma}$ of limit set of $\Gamma$ 
satisfies $\delta_{\Gamma} < 1$, 
then $F_{1}(\Gamma)$ has the following absolutely convergent product: 
$$
F_{1}(\Gamma) = \prod_{\{ \gamma \}} \prod_{m=0}^{\infty} 
\left( 1 - q_{\gamma}^{1+m} \right), 
$$
where $\{ \gamma \}$ runs over primitive conjugacy classes in $\Gamma - \{ 1 \}$, 
and $q_{\gamma}$ denotes the multiplier of $\gamma$. 

\end{itemize}
\end{it}

By a result of Krasnov [Kr] (see also [TT, (1.13)]), 
if $\Gamma$ is a Schottky group of rank $g$ 
then for $X = {\mathbb H}^{3}/\Gamma$ and $R = R_{\Gamma}$, 
$$
{\rm Vol}_{\rm R}(X) = - \frac{1}{4} S_{\rm L} - \frac{\pi}{2} \chi(R) 
= - \frac{1}{4} S_{\rm L} + \pi (g-1), 
$$
and hence one has: 
\vspace{2ex} 

{\sc Corollary 2.2}. 
\begin{it}
Let the notation be as above. 
Then 
$$
\frac{\det' \Delta_{h}}{\det {\rm Im}(\tau)} = 
c_{g} \cdot e^{(1-g)/3} 
\exp \left( \frac{{\rm Vol}_{\rm R}(X)}{3 \pi} \right) |F_{1}(\Gamma)|^{2}. 
$$
\end{it}

{\sc Theorem 2.3} (Guillarmou and Moroianu [GM, Theorem 43]). 
\begin{it} 
There exists a connection $\nabla^{\cal L}$ and a hermitian metric 
$\left\| \cdot \right\|_{\cal L}$ on the Chern-Simons line bundle 
${\cal L} = {\cal L}_{{\mathfrak S}_{g}}$ on ${\mathfrak S}_{g}$ such that 
$\lambda_{1}^{\otimes 6}$ is isomorphic to ${\cal L}_{{\mathfrak S}_{g}}^{\otimes (-1)}$ 
with respect to their connections and hermitian products  induced from those of 
$(\lambda_{1}, \nabla^{\det}, \left\| \cdot \right\|_{\rm Q})$ and 
$({\cal L}_{{\mathfrak S}_{g}}, \nabla^{\cal L}, 
\left\| \cdot \right\|_{\cal L})$ 
respectively. 
More precisely, there is an explicit isometric isomorphism 
of holomorphic hermitian line bundles given by 
$$
\left( \sqrt{c_{g} \cdot e^{1-g}} \cdot F_{1} \varphi \right)^{\otimes 6} 
\mapsto e^{-2 \pi \sqrt{-1} {\rm CS}^{PSL_{2}({\mathbb C})}}, 
$$
where $F_{1}$ is the holomorphic function on ${\mathfrak S}_{g}$, 
and $\varphi$ is the canonical section of $\lambda_{1}$ as above. 
\end{it} 
\vspace{2ex}

2.5. {\it Chern-Simons line bundle on the moduli of curves}. 
For each $0 \leq i \leq [g/2]$, 
let $t_{i}$ be a degeneration parameter of a family of 
stable (algebraic) curves of of genus $g$ such that 
the degeneration under $t_{i} = 0$ is desingularized to 
a stable curve of genus $g - 1$ if $i = 0$, 
and to $2$ stable curves of genus $i$ and $g - i$ if $i > 0$. 
\vspace{2ex}

{\sc Proposition 2.4}. 
\begin{it} 
For a family of Schottky uniformized $3$-manifolds $X$ whose boundaries $\partial X$ 
are Riemann surfaces of genus $g$ degenerating as $t_{i} \rightarrow 0$, 
one has 
$$
\left| F_{1}(\partial X)^{6} 
\exp \left( 2 \pi \sqrt{-1} {\rm CS}^{PSL_{2}({\mathbb C})} \right) \right| 
\sim |t_{i}|^{1/2} \ (t_{i} \rightarrow 0) 
$$
which means that 
$$
\lim_{t_{i} \rightarrow 0} 
\left| F_{1}(\partial X)^{6} 
\exp \left (2 \pi \sqrt{-1} {\rm CS}^{PSL_{2}({\mathbb C})} \right) \right|  
\cdot |t_{i}|^{-1/2} 
$$
exists and is not equal to $0$. 
\end{it}
\vspace{2ex} 

{\it Proof.} 
First, for each $0 \leq i \leq [g/2]$, 
we consider the associated degeneration of Schottky uniformized Riemann surfaces, 
and show that $F_{1}$ converges under this degeneration.  
Let $\Gamma = \langle \gamma_{1}, ..., \gamma_{g} \rangle$ be a Schottky group of rank $g$. 
When $i = 0$, let $\gamma_{g}(t_{0})$ be an element of $PSL_{2}({\mathbb C})$ 
which has fixed points same to those of $\gamma_{g}$ and multiplier $t_{0}$. 
Then for sufficiently small $t_{0} \neq 0$, 
$$
\Gamma(t_{0}) = \left\langle \gamma_{1}, ..., \gamma_{g-1}, \gamma_{g}(t_{0}) \right\rangle 
$$
is a Schottky group of rank $g$, and under $t_{0} \rightarrow 0$, 
the associated Riemann surface $R_{\Gamma(t_{0})}$ tends to the stable curve 
obtained from $R_{\langle \gamma_{1}, ..., \gamma_{g-1} \rangle}$ 
by identifying two fixed points of $\gamma_{g}$. 
Furthermore, the multiplier $q_{\gamma}$ of any 
$\gamma \in \Gamma(t_{0}) - \langle \gamma_{1}, ..., \gamma_{g-1} \rangle$ 
tends to $0$ as $t_{0} \rightarrow 0$, 
and hence 
$$
F_{1} \left( \Gamma(t_{0}) \right) \rightarrow 
F_{1 }\left( \langle \gamma_{1}, ..., \gamma_{g-1} \rangle \right) \sim 1 \ (t_{0} \rightarrow 0). 
$$
When $i > 0$, take two points $a, a' \in {\mathbb P}^{1}({\mathbb C})$ which are outside 
the Jordan curves $C_{\pm j}$ associated with $\gamma_{j}$ $(1 \leq j \leq g)$, 
and let $\mu$ (resp. $\mu'$) be elements of $PSL_{2}({\mathbb C})$ 
with attractive fixed point $a$ (resp. $a'$), repulsive fixed point $a'$ (resp. $a$) 
and multiplier $t_{i}$. 
Then for sufficiently small $t_{i} \neq 0$, 
$$
\Gamma(t_{i}) = \left\langle \gamma_{1}, ..., \gamma_{i}, 
\mu' \gamma_{i+1} \mu, ..., \mu' \gamma_{g} \mu \right\rangle 
$$
are Schottky groups of rank $g$, and under $t_{i} \rightarrow 0$, 
$R_{\Gamma(t_{i})}$ tends to the union of $R_{\langle \gamma_{1}, ..., \gamma_{i} \rangle}$ and 
$R_{\langle \gamma_{i+1}, ..., \gamma_{g} \rangle}$ obtained by identifying $a$ and $a'$. 
Furthermore, the multiplier $q_{\gamma}$ of any 
$\gamma \in \Gamma(t_{i}) - \langle \gamma_{1}, ..., \gamma_{i} \rangle$ 
tends to $0$ as $t_{i} \rightarrow 0$, 
and hence 
$$
F_{1}(\Gamma(t_{i})) \rightarrow F_{1}(\langle \gamma_{1}, ..., \gamma_{i} \rangle) 
\sim 1 \ (t_{i} \rightarrow 0). 
$$

Therefore, to prove the assertion, it is enough to show that for each $0 \leq i \leq [g/2]$, 
$$
\left| \exp \left( 2 \pi \sqrt{-1} {\rm CS}^{PSL_{2}({\mathbb C})} \right) \right| \sim 
|t_{i}|^{1/2} \ 
(t_{i} \rightarrow 0). 
$$
Since ${\rm CS}(g, S)$ is ${\mathbb R}$-valued, 
$$
\left| \exp \left( 2 \pi \sqrt{-1} {\rm CS}^{PSL_{2}({\mathbb C})} \right) \right| 
= e^{{\rm Vol}_{\rm R}(X)/\pi} \cdot e^{-\chi(\partial X)} 
= e^{-S_{\rm L}/(4 \pi)} \cdot e^{3(g-1)}. 
$$
If $i = 0$, then [Z1, Theorem 2.4] implies that 
$\exp \left( -S_{\rm L}/(4 \pi) \right) \sim |t_{0}|^{1/2}$. 
If $i > 0$, then Theorem 2.1 implies that 
$$
\exp \left( - S_{\rm L}/(4 \pi) \right) \sim \left\| \varphi \right\|_{\rm Q}^{-6} \ 
(t_{i} \rightarrow 0), 
$$
and by [Fr, Corollary 5.8], 
$\left\| \varphi \right\|_{\rm Q}^{-6} \sim |t_{i}|^{1/2}$ as $t_{i} \rightarrow 0$. 
This completes the proof. \ $\square$ 
\vspace{2ex} 

Let ${\cal M}_{g}$ denote the moduli stack over ${\mathbb Z}$ of 
proper smooth algebraic curves of genus $g$. 
Then the associated complex orbifold ${\cal M}_{g}({\mathbb C})$ 
becomes the moduli space of Riemann surfaces of genus $g$. 
Let $\overline{\cal M}_{g}$ denote the Deligne-Mumford compactification 
of ${\cal M}_{g}$ as the moduli stack of stable curves of genus $g$ (cf. [DM]). 
Then the complement 
$\partial {\cal M}_{g} = \overline{\cal M}_{g} - {\cal M}_{g}$ is 
the union of normal crossing divisors ${\cal D}_{i}$ 
defined as $t_{i} = 0$ $(0 \leq i \leq [g/2])$. 
For the universal stable curve 
$\pi : \overline{\cal C}_{g} \rightarrow \overline{\cal M}_{g}$, 
the Hodge line bundle $\lambda_{1} = \lambda_{g;1}$ is defined as 
$\det \left( \pi_{*} \left( \omega_{\overline{\cal C}_{g}/\overline{\cal M}_{g}} \right) \right)$, 
where $\omega_{\overline{\cal C}_{g}/\overline{\cal M}_{g}}$ denotes 
the dualizing sheaf on $\overline{\cal C}_{g}$ over $\overline{\cal M}_{g}$. 
\vspace{2ex} 

{\sc Theorem 2.5}. 
\begin{it} 
The bundle ${\cal L}_{{\mathfrak S}_{g}}^{\otimes 2}$ can be descended to 
a line bundle on ${\cal M}_{g}({\mathbb C})$ 
which we denote by ${\cal L}_{{\cal M}_{g}({\mathbb C})}$. 
Furthermore,  there exists is an isomorphism 
$$
\lambda_{g;1}^{\otimes 12} \stackrel{\sim}{\rightarrow} 
{\cal L}_{{\cal M}_{g}({\mathbb C})}^{\otimes (-1)}
$$
which is also an isometry up to a nonzero constant. 
\end{it} 
\vspace{2ex}

{\it Proof.} 
By Theorem 2.3, 
$$
\Phi = 
\left( F_{1} \varphi \right)^{12} 
\cdot \exp \left( 2 \pi \sqrt{-1} {\rm CS}^{PSL_{2}({\mathbb C})} \right)^{2} 
$$
gives a parallel section of the trivial bundle 
$\lambda_{g;1}^{\otimes 12} \otimes {\cal L}_{{\mathfrak S}_{g}}^{\otimes 2}$ 
with canonical connection. 
Since ${\mathfrak S}_{g}$ is a unramified covering 
of ${\cal M}_{g}({\mathbb C})$, 
$\Phi$ gives a local system on ${\cal M}_{g}({\mathbb C})$ 
with coefficients in ${\mathbb C}$ 
which gives rise to the monodromy representation 
$\pi_{1} \left( {\cal M}_{g}({\mathbb C}) \right) \rightarrow {\mathbb C}^{\times}$. 
Then by Proposition 2.4,  
$$
|F_{1}|^{12} \cdot 
\left| \exp \left( 2 \pi \sqrt{-1} {\rm CS}^{PSL_{2}({\mathbb C})} \right) \right|^{2} 
\sim |t_{i}| \ 
(t_{i} \rightarrow 0), 
$$
and hence a small loop in ${\cal M}_{g}({\mathbb C})$ around $t_{i} = 0$ 
has the trivial holonomy. 
Since $\pi_{1} \left( {\cal M}_{g}({\mathbb C}) \right)$ is the mapping class group of genus $g$ 
and is generated by Dehn twists, 
this monodromy representation is trivial. 
Therefore, ${\cal L}_{{\mathfrak S}_{g}}^{\otimes 2}$ can be descended to 
a hermitian holomorphic line bundle on ${\cal M}_{g}({\mathbb C})$ 
which we denote by ${\cal L}_{{\cal M}_{g}({\mathbb C})}$ 
such that $\Phi$ gives an isometric isomorphism 
$\lambda_{g;1}^{\otimes 12} \stackrel{\sim}{\rightarrow} 
{\cal L}_{{\cal M}_{g}({\mathbb C})}^{\otimes (-1)}$ up to a nonzero constant. 
\ $\square$ 
\vspace{2ex}

{\sc Corollary 2.6}. 
\begin{it} 
There exists a natural extension 
${\cal L}_{\overline{\cal M}_{g}({\mathbb C})}$ of 
${\cal L}_{{\cal M}_{g}({\mathbb C})}$ 
as a line bundle on $\overline{\cal M}_{g}({\mathbb C})$. 
Furthermore, 
$$
\lambda_{g;1}^{\otimes 12} \cong 
{\cal L}_{\overline{\cal M}_{g}({\mathbb C})}^{\otimes (-1)} \otimes 
{\cal O}_{\overline{\cal M}_{g}({\mathbb C})} 
\left( \partial {\cal M}_{g} \right) 
$$
which associates $(F_{1} \varphi)^{12}$ with 
$\exp \left( 2 \pi \sqrt{-1} {\rm CS}^{PSL_{2}({\mathbb C})} \right)^{-2}$ 
up to a nonzero constant. 
\end{it}
\vspace{2ex} 

{\it Proof.} 
This follows immediately from Theorem 2.5 and its proof. 
\ $\square$. 
\vspace{2ex}

\begin{center}
3. {\bf Arithmetic Riemann-Roch theorem} 
\end{center}

3.1. {\it Riemann-Roch isomorphism}. 
Fix integers $g, n \geq 0$ such that $2g - 2 + n > 0$, 
and let $\overline{\cal M}_{g,n}$ be the moduli stack over ${\mathbb Z}$ 
of stable $n$-pointed curves of genus $g$ (cf. [K]). 
Then by definition, 
there exist the universal curve 
$\pi : \overline{\cal C}_{g,n} \rightarrow \overline{\cal M}_{g,n}$ 
and the universal sections 
$\sigma_{j} : \overline{\cal M}_{g,n} \rightarrow \overline{\cal C}_{g,n}$ 
$(1 \leq j \leq n)$. 
Denote by ${\cal M}_{g,n}$ the open substack of $\overline{\cal M}_{g,n}$ 
classifying proper smooth $n$-pointed curves of genus $g$. 
Then from the relative dualizing sheaf $\omega_{\overline{\cal C}_{g,n}/\overline{\cal M}_{g,n}}$ 
and the complement $\partial {\cal M}_{g,n} = \overline{\cal M}_{g,n} - {\cal M}_{g,n}$ 
of ${\cal M}_{g,n}$, 
one has the following line bundles on $\overline{\cal M}_{g,n}$: 
\begin{eqnarray*}
\lambda_{g,n;k} 
& : = & 
\det R \pi_{*} \left( 
\omega^{k}_{\overline{\cal C}_{g,n}/\overline{\cal M}_{g,n}} 
\left( k \sum_{j} \sigma_{j} \right) \right), 
\\ 
\psi_{g,n} 
& : = & 
\bigotimes_{j=1}^{n} \psi_{g,n}^{(j)}; \ 
\psi_{g,n}^{(j)} : = \sigma^{*}_{j} 
\left( \omega_{\overline{\cal C}_{g,n}/\overline{\cal M}_{g,n}} \right), 
\\ 
\delta_{g,n} 
& : = & 
{\cal O}_{\overline{\cal M}_{g,n}} \left( \partial {\cal M}_{g,n} \right). 
\end{eqnarray*} 
Furthermore, 
let $\kappa_{g,n}$ be the line bundle on $\overline{\cal M}_{g,n}$ defined as 
the following Deligne's pairing: 
$$
\kappa_{g,n} = 
\left\langle 
\omega_{\overline{\cal C}_{g,n}/\overline{\cal M}_{g,n}} 
\left( \sum_{j} \sigma_{j} \right), 
\omega_{\overline{\cal C}_{g,n}/\overline{\cal M}_{g,n}} 
\left( \sum_{j} \sigma_{j} \right) 
\right\rangle. 
$$
Using the $k$th residue map ${\rm Res}^{k}_{\sigma_{j}}$ given by 
${\rm Res}^{k}_{\sigma_{j}} \left( \eta \left( dz_{j}/z_{j} \right)^{k} 
\right) = \eta$, 
where $z_{j}$ denotes a local coordinate on 
$\overline{\cal C}_{g,n}$ around $\sigma_{j}$, 
one has 
$$
\lambda_{g,n;k} \cong 
\det \pi_{*} \left( 
\omega^{k}_{\overline{\cal C}_{g,n}/\overline{\cal M}_{g,n}} 
\left( (k - 1) \sum_{j} \sigma_{j} \right) \right). 
$$
In particular, 
$\lambda_{g,n;1} \cong \det \pi_{*} 
\left( \omega_{\overline{\cal C}_{g,n}/\overline{\cal M}_{g,n}} \right)$. 
When $n = 0$, 
we delete $n$ from these notations, 
and put $\psi_{g,0} = \psi_{g} = {\cal O}_{\overline{\cal M}_{g}}$. 
Furthermore, put $d_{k} = 6 k^{2} - 6 k + 1$ and 
$$
a(g, n) = (2g - 2 + n) \left( -12 \zeta'_{\mathbb Q}(-1) + \frac{1}{2} \right), 
$$
where $\zeta'_{\mathbb Q}(-1)$ denotes the derivative of Riemann's zeta function  $\zeta_{\mathbb Q}$ at $-1$. 
\vspace{2ex} 

{\sc Theorem 3.1} (arithmetic Riemann-Roch theorem [D2, Fr, GS, W]). 
\begin{it}
There exists a unique (up to a sign) isomorphism 
$$
\lambda_{g,n;k}^{\otimes 12} \otimes \psi_{g,n} \otimes 
\delta_{g,n}^{\otimes (-1)} \cong \kappa_{g,n}^{\otimes d_{k}} \cdot e^{a(g,n)} 
$$
between the line bundles over $\overline{\cal M}_{g,n}$ 
which is an isometry between the line bundles over ${\cal M}_{g,n}({\mathbb C})$ 
for these hermitian structure. 
\end{it}
\vspace{2ex}

3.2. {\it Mumford isomorphism}. 
\vspace{2ex}

{\sc Theorem 3.2}. 
\begin{it} 
There exists a unique (up to a sign) isomorphism 
$$
\mu_{g,n;k} : \lambda_{g,n;k} \stackrel{\sim}{\rightarrow} 
\lambda_{g,n;1}^{\otimes d_{k}} \otimes 
\left( \psi_{g,n} \otimes \delta_{g,n}^{\otimes (-1)} 
\right)^{\otimes (k^{2}-k)/2} \cdot \exp \left( (k-k^{2}) a(g, n)/2 \right). 
$$
between the line bundles over $\overline{\cal M}_{g,n}$ 
which is an isometry between the line bundles over ${\cal M}_{g,n}({\mathbb C})$ 
for these hermitian structure. 
We call $\mu_{g,n;k}$ the Mumford isomorphism (cf. [Mu2]). 
\end{it} 
\vspace{2ex}

{\it Proof.} 
The uniqueness follows from the properness of $\overline{\cal M}_{g,n}$ 
over ${\mathbb Z}$. 
The existence is shown by Mumford [Mu2] when $n = 0$, 
namely, there exists a canonical isomorphism 
$$
\lambda_{g+n;k} \cong \lambda_{g+n;1}^{\otimes d_{k}} \otimes 
\delta_{g+n}^{\otimes (k - k^{2})/2}. 
$$
Let 
${\rm cl}_{1} : 
\overline{\cal M}_{g,n} \times \overline{\cal M}_{1,1}^{\times n} 
\rightarrow \overline{\cal M}_{g+n}$ 
be the clutching morphism 
given by identifying the $j$th section $\sigma_{j}$ 
over $\overline{\cal M}_{g,n}$ and the unique section of the $j$th component 
of $\overline{\cal M}_{1,1}^{\times n}$ $(1 \leq j \leq n)$. 
Then by [Fr, Corollary 3.4], 
\begin{eqnarray*} 
{\rm cl}_{1}^{*} \left( \lambda_{g+n;k} \right) 
& \cong & 
\lambda_{g,n;k} \boxtimes \lambda_{1,1;k}^{\boxtimes n}, 
\\ 
{\rm cl}_{1}^{*} \left( \delta_{g+n} \right) 
& \cong & 
\left( \delta_{g,n} \otimes \psi_{g,n}^{\otimes (-1)} \right) \boxtimes 
\left( \delta_{1,1} \otimes \psi_{1,1}^{\otimes (-1)} \right)^{\boxtimes n}, 
\end{eqnarray*} 
and there exists an isomorphism between 
$\lambda_{g,n;k} \boxtimes \lambda_{1,1;k}^{\boxtimes n}$ and 
$$
\left( \lambda_{g,n;1}^{\otimes d_{k}} \otimes \left( 
\delta_{g,n} \otimes \psi_{g,n}^{\otimes (-1)} \right)^{\otimes (k - k^{2})/2} 
\right) \boxtimes 
\left( \lambda_{1,1;1}^{\otimes d_{k}} \otimes \left( 
\delta_{1,1} \otimes \psi_{1,1}^{\otimes (-1)} \right)^{\otimes (k - k^{2})/2} 
\right)^{\boxtimes n}. 
$$
Denote by 
$\Delta : {\rm Spec}({\mathbb Z}) \rightarrow 
\overline{\cal M}_{1,1}^{\times n}$ 
the section corresponding to the $n$ copies of the $1$-pointed stable curve 
obtained from ${\mathbb P}_{\mathbb Z}^{1}$ by identifying $0$ and $\infty$ 
with marked point $1$. 
Then taking the pullback of the above isomorphism by 
${\rm id} \times \Delta : \overline{\cal M}_{g,n} \rightarrow 
\overline{\cal M}_{g,n} \times \overline{\cal M}_{1,1}^{\times n}$, 
we have the required isomorphism. 
\ $\square$ 
\vspace{2ex} 

{\sc Proposition 3.3}. 
\begin{it} 
Let 
${\rm cl}_{2} : \overline{\cal M}_{g,n+2} \rightarrow 
\overline{\cal M}_{g+1,n}$ 
be the clutching morphism given by identifying the sections 
$\sigma_{n+1}$ and $\sigma_{n+2}$. 
Then there exist the following canonical isomorphisms: 
\begin{eqnarray*}
{\rm cl}_{2}^{*} \left( \lambda_{g+1,n;k} \right) & \cong & 
\lambda_{g,n+2;k}, 
\\ 
{\rm cl}_{2}^{*} \left( \psi_{g+1,n} \right) & \cong & 
\psi_{g,n+2} \otimes \left( \psi_{g,n+2}^{(n+1)} \right)^{\otimes(-1)} 
\otimes \left( \psi_{g,n+2}^{(n+2)} \right)^{\otimes(-1)}, 
\\ 
{\rm cl}_{2}^{*} \left( \delta_{g+1,n} \right) & \cong & 
\delta_{g,n+2} \otimes \left( \psi_{g,n+2}^{(n+1)} \right)^{\otimes(-1)} 
\otimes \left( \psi_{g,n+2}^{(n+2)} \right)^{\otimes(-1)}, 
\end{eqnarray*}
and under these isomorphisms,  
${\rm cl}_{2}^{*} \left( \mu_{g+1,n;k} \right) = \pm \mu_{g,n+2;k}$. 
\end{it}
\vspace{2ex}

{\it Proof.} 
Let $(\pi : C \rightarrow S; \sigma_{1},..., \sigma_{n+2})$ 
be an $(n+2)$-pointed stable curve of genus $g$, 
and $(\pi' : C' \rightarrow S; \sigma_{1},..., \sigma_{n})$ 
be the $n$-pointed stable curve of $g + 1$ obtained from $C$ 
by identifying $\sigma_{n+1}$ and $\sigma_{n+2}$. 
Then as is shown in [K, Section 1], 
$\pi'_{*} \left( \omega^{k}_{C'/S} \left( k \sum_{j=1}^{n} \sigma_{j} 
\right) \right)$ 
is isomorphic to ${\rm Ker}(\rho)$, 
where  
$\rho : \pi_{*} \left( \omega^{k}_{C/S} \left( k \sum_{j=1}^{n+2} \sigma_{j} 
\right) \right) \rightarrow {\cal O}_{S}$ 
is given by 
$$
\rho(\eta) = {\rm Res}^{k}_{\sigma_{n+1}} (\eta) - 
(-1)^{k} {\rm Res}^{k}_{\sigma_{n+2}} (\eta). 
$$
Therefore, we have the first isomorphism. 
The second and third ones are shown in [K, Theorem 4.3], 
and hence the last identity follows from 
the uniqueness of the Mumford isomorphism. 
\ $\square$ 
\vspace{2ex}

\begin{center}
4. {\bf Arithmetic Schottky uniformization} 
\end{center}

4.1. {\it Tate curve}. 
According to [Si, T], we review the Tate curve over the ring ${\mathbb Z}[[q]]$ of integral power series of $q$ which gives the universal elliptic curve over the ring  ${\mathbb Z}((q)) = {\mathbb Z}[[q]][1/q]$ of integral Laurent power series of $q$. 
Put 
$$
s_{k} = \sum_{n=1}^{\infty} \frac{n^{k} q^{n}}{1 - q^{n}}, \ \ 
a_{4}(q) = - 5 s_{3}(q), \ \ 
a_{6}(q) = - \frac{5 s_{3}(q) + 7 s_{5}(q)}{12} 
$$
which are seen to be in ${\mathbb Z}[[q]]$. 
Then the Tate curve ${\cal E}_{q}$ is defined as 
$$
y^{2} + xy = x^{3} + a_{4}(q) x + a_{6}(q), 
$$
and is formally represented as the quotient space ${\mathbb G}_{m} / \langle q \rangle$ 
with the origin $o$ given by $1 \in {\mathbb G}_{m}$. 
Therefore, 
$dz/z$ $(z \in {\mathbb G}_{m})$ is a regular $1$-form on ${\cal E}_{q}$, 
and 
\begin{eqnarray*}
X(z) & = & 
\sum_{n \in {\mathbb Z}} \frac{q^{n} z}{(1 - q^{n} z)^{2}} - 2 s_{1}(q) 
\\ & = & 
\frac{z}{(1 - z)^{2}} + \sum_{n=1}^{\infty} 
\left( \frac{q^{n} z}{(1 - q^{n} z)^{2}} + 
\frac{q^{n} z^{-1}}{(1 - q^{n} z^{-1})^{2}} - 
2 \frac{q^{n}}{(1 - q^{n})^{2}} \right), 
\\ 
Y(z) & = & 
\sum_{n \in {\mathbb Z}} \frac{(q^{n} z)^{2}}{(1 - q^{n} z)^{3}} + s_{1}(q) 
\\ & = & 
\frac{z^{2}}{(1 - z)^{3}} + \sum_{n=1}^{\infty} 
\left( \frac{(q^{n} z)^{2}}{(1 - q^{n} z)^{3}} + 
\frac{q^{n} z^{-1}}{(1 - q^{n} z^{-1})^{3}} + 
\frac{q^{n}}{(1 - q^{n})^{2}} \right) 
\end{eqnarray*}
are meromorphic functions on ${\cal E}_{q}$ which have only one pole at $o$ 
of order $2$, $3$ respectively. 
Hence for each positive integer $k$, 
$H^{0} \left( {\cal E}_{q}, \omega^{k}_{{\cal E}_{q}/{\mathbb Z}[[q]]} (ko) \right)$ 
is a free ${\mathbb Z}[[q]]$-module of rank $k$ generated by 
$(dz/z)^{k}$ and $f_{i}(z) (dz/z)^{k}$ $(2 \leq i \leq k)$, 
where $f_{i}(z)$ is a meromorphic function on ${\cal E}_{q}$ 
with only one pole at $o$ of order $i$. 
\vspace{2ex}

We consider the Tate curve with marked points. 
For variables $t_{1},..., t_{h}$, 
put 
$$
R = {\mathbb Z} \left[ t_{k}, \frac{1}{t_{k}}, \frac{1}{t_{k}-1}, 
\frac{1}{t_{l} - t_{m}} \ (1 \leq k, l, m \leq h, \ l \neq m) \right]. 
$$
Then each $t_{i}$ gives a point on ${\cal E}_{q} \otimes R[[q]]$ 
which we denote by the same symbol. 
When $t_{1} =1$, the corresponding point becomes the origin $o$. 
\vspace{2ex} 

{\sc Proposition 4.1}. 
\begin{it} 
The $R[[q]]$-module  
$$
V = H^{0} \left( {\cal E}_{q} \otimes R[[q]], 
\Omega^{k}_{{\cal E}_{q} \otimes R[[q]]/R[[q]]} 
\left( k \sum_{j=1}^{h} t_{j} \right) \right)
$$
is free of rank $hk$, and its basis consists of certain products of 
$X(z/t_{j})/t_{j}^{2}$, $Y(z/t_{j})/t_{j}^{3}$ times $(dz/z)^{k}$. 
\end{it}
\vspace{2ex}

{\it Proof.} 
This follows from that ${\rm Res}_{z=0} \left( z^{k-1} (dz/z^{k}) \right) = 1$, 
and that $X(z/t_{j})/t_{j}^{2}$, $Y(z/t_{j})/t_{j}^{3}$ have Laurent expansion 
at $z = 0$ with coefficients in $R((q))$. 
\ $\square$ 
\vspace{2ex}

4.2. {\it Arithmetic Schottky uniformization}. 
The arithmetic Schottky uniformization theory [I2] gives 
a higher genus version of the Tate curve, 
and its $1$-forms and periods. 
We review this theory for the special case concerned with 
universal deformations of irreducible degenerate curves. 

Denote by $\Delta$ the graph with one vertex and $g$ loops. 
Let $x_{\pm 1},..., x_{\pm g}$, $y_{1},..., y_{g}$ be variables, 
and put 
\begin{eqnarray*}
A_{g} & = & {\mathbb Z} \left[ x_{k}, \frac{1}{x_{l}-x_{m}} \ 
\left( k, l, m \in \{ \pm 1,..., \pm g \}, \ l \neq m \right) \right], 
\\
A_{\Delta} & = & A [[ y_{1},... y_{g} ]], 
\\
B_{\Delta} & = & A_{\Delta} \left[ 1/y_{i} \ (1 \leq i \leq g) \right]. 
\end{eqnarray*}
Then it is shown in [I2, Section 3] that 
there exists a stable curve $C_{\Delta}$ of genus $g$ over $A_{\Delta}$ 
which satisfies the followings: 
\begin{itemize}

\item 
$C_{\Delta}$ is a universal deformation of the universal degenerate curve 
with dual graph $\Delta$ which is obtained from ${\mathbb P}^{1}_{A}$ 
by identifying $x_{i}$ and $x_{-i}$ $(1 \leq i \leq g)$. 
The ideal of $A_{\Delta}$ generated by $y_{1},..., y_{g}$ 
corresponds to the boundary $\partial {\cal M}_{g} = \overline{\cal M}_{g} - {\cal M}_{g}$ 
of $\overline{\cal M}_{g}$ via the morphism 
${\rm Spec}(A_{\Delta}) \rightarrow \overline{\cal M}_{g}$ 
associated with $C_{\Delta}$. 

\item 
$C_{\Delta}$ is smooth over $B_{\Delta}$, 
and is Mumford uniformized (cf. [Mu1]) by the subgroup $\Gamma_{\Delta}$ of 
$PGL_{2} \left( B_{\Delta} \right)$ 
with $g$ generators 
$$
\phi_{i} = 
\left( \begin{array}{cc} x_{i} & x_{-i} \\ 1 & 1 \end{array} \right) 
\left( \begin{array}{cc} y_{i} & 0      \\ 0 & 1 \end{array} \right) 
\left( \begin{array}{cc} x_{i} & x_{-i} \\ 1 & 1 \end{array} \right)^{-1} 
\ {\rm mod} \left( B_{\Delta}^{\times} \right) 
\ (1 \leq i \leq g). 
$$
Furthermore, 
$C_{\Delta}$ has the following universality: 
for a complete integrally closed noetherian local ring $R$ 
with quotient field $K$ 
and a Mumford curve $C$ over $K$ such that $\Delta$ is the dual graph of 
its degenerate reduction, 
there is a ring homomorphism $A_{\Delta} \rightarrow R$ 
which gives rise to $C_{\Delta} \otimes_{A_{\Delta}} K \cong C$. 

\item 
Let $\Gamma = \langle \gamma_{1},..., \gamma_{g} \rangle$ be a Schottky group of rank $g$, 
where each $\gamma_{i}$ has the attractive (resp. repulsive) fixed points $a_{i}$ (resp. $a_{-i}$), 
and it has the multiplier $q_{i}$. 
Then substituting $a_{\pm i}$ to $x_{\pm i}$ and $q_{i}$ to $y_{i}$ $(1 \leq i \leq g)$, 
$C_{\Delta}$ becomes the Riemann surface $R_{\Gamma}$ uniformized by $\Gamma$ 
if $|q_{i}|$ are sufficiently small. 

\end{itemize}

Actually, $C_{\Delta}$ is constructed in [I2] as the quotient of 
a certain subspace of ${\mathbb P}^{1}_{B_{\Delta}}$ by the action of 
$\Gamma$ using the theory of formal schemes. 
Furthermore, as is shown in [MD] and [I1, Section 3], 
there exists a basis of sections  
\begin{eqnarray*} 
\omega_{i} 
& = & 
\sum_{\phi \in \Gamma_{\Delta} / \langle \phi_{i} \rangle} 
\left( \frac{1}{z - \phi(x_{i})} - \frac{1}{z - \phi(x_{-i})} \right) 
\\
& = & 
\left( \frac{1}{z - x_{i}} - \frac{1}{z - x_{-i}} \right) + 
\sum_{\phi \in \Gamma_{\Delta} / \langle \phi_{i} \rangle - \{ 1 \}} 
\left( \frac{1}{z - \phi(x_{i})} - \frac{1}{z - \phi(x_{-i})} \right) + 
\cdots 
\\
& \in & 
A_{\Delta} \left[ \prod_{k=1}^{g} \frac{1}{(z - x_{k})(z - x_{-k})} \right] \ (1 \leq i \leq g) 
\end{eqnarray*}
of the dualizing sheaf $\omega_{C_{\Delta}/A_{\Delta}}$ on $C_{\Delta}$. 
These $\omega_{i}$ are called {\it normalized} since they give rise to 
holomorphic $1$-forms on the above $R_{\Gamma}$ 
such that $\oint_{C_{j}} \omega_{i} = 2 \pi \sqrt{-1} \delta_{ij}$, 
where $C_{j} \subset {\mathbb P}^{1}({\mathbb C})$ are given in 2.2. 
Put $\omega = \omega_{1} \wedge \cdots \wedge \omega_{g}$. 

Let $\Delta_{g-1}$ be the graph with one vertex and $(g-1)$ loops, 
and put 
$$
A_{g-1} = {\mathbb Z} \left[ x_{k}, \frac{1}{x_{l}-x_{m}} \ 
\left( k, l, m \in \{ \pm 1,..., \pm (g-1) \}, \ l \neq m \right) \right]. 
$$
Then one has the generalized Tate curve $C_{\Delta_{g-1}}$ of genus $g-1$ 
over $A_{g-1}[[ y_{1},..., y_{g-1}]]$. 
\vspace{2ex} 

{\sc Proposition 4.2}. 
\begin{it} 
\begin{itemize}

\item[{\rm (1)}] 
The stable curve $C_{\Delta}|_{y_{g} = 0}$ is obtained from 
$C_{\Delta_{g-1}} \otimes_{A_{g-1}} A_{g}$ by identifying $x_{g} = x_{-g}$. 

\item[{\rm (2)}] 
The set $\left\{ \omega_{i}|_{y_{g} =0} \right\}_{1 \leq i \leq g}$ gives a basis of sections of 
the dualizing sheaf on $ C_{\Delta}|_{y_{g}=0}$ over $A_{g}[[y_{1},..., y_{g-1}]]$ which are 
normalized in the sense that $\omega_{i}|_{y_{g} = 0}$ $(i < g)$ are normalized $1$-forms 
on $C_{\Delta_{g-1}}$ over $A_{g-1}[[y_{1},..., y_{g-1}]]$ and 
${\rm Res}_{z = x_{g}} \left( \omega_{g}|_{y_{g} = 0} \right) = 1$. 

\end{itemize}
\end{it} 

{\it Proof.} 
The assertion (1) follows from the universality of $C_{\Delta}$ and $C_{\Delta_{g-1}}$, 
and (2) follows from the above universal expression of $\omega_{i}$. 
\ $\square$ 
\vspace{2ex}

In the case when we consider the Schottky space ${\mathfrak S}_{g}$ 
of degree $g > 1$ as the moduli space of normalized Schottky groups, 
we assume that the above $\phi_{1},..., \phi_{g}$ are {\it normalized} 
by considering $x_{1}, x_{-1}$ as $0, \infty$ respectively, namely, 
$$
\phi_{1} = 
\left( \begin{array}{cc} 1 & 0 \\ 0 & y_{1} \end{array} \right) 
\ {\rm mod} \left( B_{\Delta}^{\times} \right), 
$$
and by putting $x_{2} = 1$. 
Then as is shown in [I2, 1.1], the corresponding generalized Tate curve $C_{\Delta}$ 
is defined over $A'_{\Delta} = A'_{g} [[y_{1},..., y_{g}]]$, 
where $A'_{g}$ is obtained from $A_{g}$ by deleting $x_{-1}$ and putting $x_{1} = 0$, $x_{2} = 1$. 
Therefore, one has the associated morphism 
${\rm Spec} \left( A'_{\Delta} \right) \rightarrow \overline{\cal M}_{g}$. 
\vspace{2ex}

{\sc Proposition 4.3}. 
\begin{it} 
\begin{itemize}

\item[{\rm (1)}] 
The infinite products 
$$
\prod_{\{ \gamma \}} \prod_{m=0}^{\infty} \left( 1 - q_{\gamma}^{1+m} \right), 
\ \ 
\left( 1 - q_{\gamma_{1}} \right)^{2} \cdots 
\left( 1 - q_{\gamma_{1}}^{k-1} \right)^{2} 
\left( 1 - q_{\gamma_{2}}^{k-1} \right) 
\prod_{\{ \gamma \}} \prod_{m=0}^{\infty} \left( 1 - q_{\gamma}^{k+m} \right) 
$$
have universal expression as elements of $A'_{\Delta}$. 

\item[{\rm (2)}] 
Under $y_{2} = \cdots = y_{g} = 0$, 
these elements of $A'_{\Delta}$ becomes 
$\prod_{m=0}^{\infty} \left( 1 - y_{1}^{1+m} \right)^{2}$ 
which is primitive, 
i.e., not congruent to $0$ modulo any rational prime. 

\end{itemize}
\end{it}

{\it Proof.} 
Let $(\Gamma; \gamma_{1},..., \gamma_{g})$ be a normalized Schottky group, 
and for $i = 1,..., g$, put $\gamma_{-i} = \gamma_{i}^{-1}$. 
Then by Proposition 1.3 of [I2] and its proof, 
if $\gamma \in \Gamma - \{ 1 \}$ has the reduced expression 
$\gamma_{\sigma(1)} \cdots \gamma_{\sigma(l)}$ $(\sigma(i) \in \{ \pm 1,..., \pm g \})$ 
such that $\sigma(1) \neq - \sigma(l)$, 
then its multiplier $q_{\gamma}$ has universal expression as an element 
of $A'_{\Delta}$ divisible by $y_{\sigma(1)} \cdots y_{\sigma(l)}$. 
Therefore, the assertion (1) holds. 
Since 
$$
\left( 1 - y_{1} \right)^{2} \cdots \left( 1 - y_{1}^{k-1} \right)^{2} 
\prod_{m=0}^{\infty} \left( 1 - y_{1}^{k+m} \right)^{2}  
= \prod_{m=0}^{\infty} \left( 1 - y_{1}^{1+m} \right)^{2}, 
$$
the assertion (2) holds. 
\ $\square$ 
\vspace{2ex}

\begin{center}
5. {\bf Explicit Riemann-Roch isomorphism} 
\end{center}

5.1. {\it Arithmeticity of Chern-Simons invariant}. 
\vspace{2ex}

{\sc Theorem 5.1}. 
\begin{it} 
\begin{itemize}

\item[{\rm (1)}] 
There exists an isomorphism 
${\cal L}_{\overline{\cal M}_{g}({\mathbb C})} \cong \kappa_{g}^{\otimes (-1)}$ 
between the line bundles over $\overline{\cal M}_{g}({\mathbb C})$ 
which is isometric over ${\cal M}_{g}({\mathbb C})$. 

\item[{\rm (2)}] 
There exists a unique (up to a sign) primitive element of $A_{\Delta}$ 
which gives a universal expression of 
$$
\exp \left( 4 \pi \sqrt{-1} {\rm CS}^{PSL_{2}({\mathbb C})} \right) 
$$
times a certain nonzero constant. 

\end{itemize}
\end{it}

{\it Proof.} 
The assertion (1) follows from Corollary 2.6 and Theorem 3.1 
(in the case when $n = 0$, $k = 1$), 
and then we prove the assertion (2). 
By the arithmetic Schottky uniformization theory, 
$\omega$ gives a trivialization of $\lambda_{g;1}$ over 
${\rm Spec} \left( A_{\Delta} \right)$. 
Since $\omega = \left( 2 \pi \sqrt{-1} \right)^{g} \varphi$, 
by the description of the isomorphism in Theorem 2.5 and Corollary 2.6, 
a certain multiple of the holomorphic function 
$$
F_{1}^{12} \cdot 
\exp \left( 4 \pi \sqrt{-1} {\rm CS}^{PSL_{2}({\mathbb C})} \right) 
$$ 
on  ${\mathfrak S}_{g}$ gives a trivialization over ${\rm Spec} \left( A_{\Delta} \right)$ 
of 
$$
{\cal L}_{\overline{\cal M}_{g}({\mathbb C})} \otimes 
{\cal O}_{\overline{\cal M}_{g}} \left( -\partial {\cal M}_{g} \right) 
\cong 
\kappa_{g}^{\otimes (-1)} \otimes 
{\cal O}_{\overline{\cal M}_{g}} \left( -\partial {\cal M}_{g} \right)  
$$
for its  ${\mathbb Z}$-structure induced from that on $\kappa_{g}$. 
As is shown in Proposition 4.3, 
$F_{1}$ is given as an invertible element of $A_{\Delta}$, 
and hence the assertion (2) holds. 
\ $\square$ 
\vspace{2ex}

{\it Definition}. 
We denote by ${\rm ACS}_{g}$ the element of $A_{\Delta}$ given in Theorem 5.1 (2), 
and call it the {\it arithmetic universal Chern-Simons invariant}. 
\vspace{2ex} 

5.2. {\it Holomorphic factorization formula}. 
Assume that $g > 1$, and take an integer $k > 1$. 
Let $\left( \Gamma; \gamma_{1},..., \gamma_{g} \right)$ be 
a marked normalized Schottky group, 
and ${\mathbb C}[z]_{2k-2}$ be the ${\mathbb C}$-vector space of polynomials 
$f = f(z)$ of $z$ with degree $\leq 2k-2$ on which $\Gamma$ acts as 
$$
\gamma(f)(z) = f \left( \gamma(z) \right) \cdot \gamma'(z)^{1-k} \ 
\left( \gamma \in \Gamma, \ f \in {\mathbb C}[z]_{2k-2} \right). 
$$
Take $\xi_{1,k-1}$, $\xi_{2,1},..., \xi_{2,2k-2}$, 
$\xi_{i,0},..., \xi_{i,2k-2}$ $(3 \leq i \leq g)$ 
as elements of the Eichler cohomology group 
$H^{1} \left( \Gamma, {\mathbb C}[z]_{2k-2} \right)$ 
of $\Gamma$ which are uniquely determined by the condition: 
$$
\xi_{i,j}(\gamma_{l}) = \left\{ \begin{array}{ll} 
\delta_{2l} (z-1)^{j} & (i = 2), 
\\ 
\delta_{il} z^{j} & (i \neq 2) 
\end{array} \right. 
$$
for $1 \leq l \leq g$. 
Then it is shown in [MT, Section 4] that 
$$
\Psi_{g;k} (\psi, \xi) := 
\frac{1}{2 \pi \sqrt{-1}} \sum_{i=1}^{g} 
\oint_{C_{i}} \psi \cdot \xi(\gamma_{i}) dz 
$$ 
for $\psi \in H^{0} \left( R_{\Gamma}, \Omega^{k}_{R_{\Gamma}} \right)$, 
$\xi \in H^{1} \left( \Gamma, {\mathbb C}[z]_{2k-2} \right)$ 
is a non-degenerate pairing on 
$$
H^{0} \left( R_{\Gamma}, \Omega^{k}_{R_{\Gamma}} \right) \times 
H^{1} \left( \Gamma, {\mathbb C}[z]_{2k-2} \right). 
$$
Then there exists a basis 
$$ 
\left\{ \psi_{1,k-1}, \ \psi_{2,1},..., \psi_{2,2k-2}, \ \psi_{i,0},..., \psi_{i,2k-2} \ 
(3 \leq i \leq g) \right\} 
$$ 
which we call the {\it normalized basis} of 
$H^{0} \left( R_{\Gamma}, \Omega^{k}_{R_{\Gamma}} \right)$ such that 
$\Psi_{g;k} \left( \psi_{i,j}, \xi_{l,m} \right) = \delta_{il} \cdot \delta_{jm}$, 
where $\delta_{ij}$ denotes Kronecker's delta. 
\vspace{2ex} 

{\it Remark.} 
Since $-\pi \cdot \Psi_{g;k}$ is the pairing given in [MT, (4.1)], 
$$ 
\left\{ -\frac{\psi_{1,k-1}}{\pi}, \ -\frac{\psi_{2,1}}{\pi},..., -\frac{\psi_{2,2k-2}}{\pi}, \ 
-\frac{\psi_{i,0}}{\pi},..., -\frac{\psi_{i,2k-2}}{\pi} \ (3 \leq i \leq g) \right\} 
$$ 
is the {\it natural basis for $n$-differentials} defined in [MT]. 
\vspace{2ex} 

In what follows, put 
$$
\left\{ \omega^{(k)}_{1},..., \omega^{(k)}_{(2k-1)(g-1)} \right\} = 
\left\{ \psi_{1,k-1}, \ \psi_{2,1},..., \psi_{2,2k-2}, \ \psi_{i,0},..., \psi_{i,2k-2} \ 
(3 \leq i \leq g) \right\}, 
$$
and $\omega^{(k)} =  \omega^{(k)}_{1} \wedge \cdots \wedge \omega^{(k)}_{(2k-1)(g-1)}$. 
\vspace{2ex} 

{\sc Theorem 5.2} (McIntyre-Takhtajan [MT, Theorem 2]). 
\begin{it} 
Assume that $k > 1$. 
\begin{itemize}

\item[{\rm (1)}] 
There exists a holomorphic function $F_{k}$ on ${\mathfrak S}_{g}$ 
which gives an isometry between $\lambda_{g;k}$ on ${\cal M}_{g}({\mathbb C})$ 
with Quillen metric and the holomorphic line bundle 
on ${\cal M}_{g}({\mathbb C})$ 
determined by the hermitian metric $\exp \left( S_{\rm L}/12 \pi \right)^{d_{k}}$. 
More precisely, 
there exists a positive real number $c_{g;k}$ depending only on $g$ and $k$ 
such that 
$$
\exp \left( \frac{S_{\rm L}}{12 \pi} \right)^{d_{k}} = 
c_{g;k} \left| F_{k} \right|^{2} \left\| \omega^{(k)} \right\|_{\rm Q}^{2}, 
$$
where $\| * \|_{\rm Q}$ denotes the Quillen metric. 

\item[{\rm (2)}] 
On the whole Schottky space ${\mathfrak S}_{g}$ classifying 
marked normalized Schottky groups 
$\left( \Gamma; \gamma_{1},..., \gamma_{g} \right)$, 
$F_{k}$ is given by the absolutely convergent infinite product 
$$
\left( 1 - q_{\gamma_{1}} \right)^{2} \cdots 
\left( 1 - q_{\gamma_{1}}^{k-1} \right)^{2} 
\left( 1 - q_{\gamma_{2}}^{k-1} \right) 
\prod_{\{ \gamma \}} \prod_{m=0}^{\infty} \left( 1 - q_{\gamma}^{k+m} \right), 
$$
where $\{ \gamma \}$ runs over primitive conjugacy classes in 
$\Gamma - \{ 1 \}$. 

\end{itemize}
\end{it}

{\it Remark.} 
As is seen in Proposition 4.3, 
$F_{k}$ has a universal expression as an element of $A'_{\Delta}$ 
which we denote by the same symbol. 
\vspace{2ex}

{\sc Proposition 5.3}. 
\begin{it} 
Denote by $\mu_{g;k}$ the Mumford isomorphism $\mu_{g,0;k}$ given in Theorem 3.2. 
Then there exists a nonzero constant $c(g;k)$ depending only on $g, k > 1$ such that 
\end{it}
$$
c(g;k) \cdot \mu_{g;k} \left( \omega^{(k)} \right) 
= \frac{F_{1}^{d_{k}}}{F_{k}} \omega^{\otimes d_{k}}. 
$$

{\it Proof.} \ 
Let $a(g)$ denote the Deligne constant 
$(1-g)\left( 24 \zeta'_{\mathbb Q}(-1) -1 \right)$.  
Then by Theorem 3.2, $\mu_{g;k}$ gives rise to an isometry 
$$
\lambda_{g;k} \cong \lambda_{g;1}^{\otimes d_{k}} \cdot \exp \left( (k - k^{2}) a(g)/2 \right) 
$$
between the metrized tautological line bundles with Quillen metric on ${\cal M}_{g}$. 
Therefore, by the formulas of Zograf and of McIntyre-Takhtajan, 
there exists a holomorphic function $c(g;k)$ on the Schottky space ${\mathfrak S}_{g}$ 
satisfying the above formula such that $\left| c(g;k) \right|$ is a constant function. 
Since ${\mathfrak S}_{g}$ is a connected complex manifold, 
$c(g;k)$ is also a constant function. 
\ $\square$ 
\vspace{2ex}
 
Let $(\Gamma; \gamma_{1},..., \gamma_{g})$ 
be a marked Schottky group of rank $g > 1$, 
and $R_{\Gamma} = \Omega_{\Gamma}/\Gamma$ be the Riemann surface 
uniformized by $\Gamma$ with $n$-marked points 
given by $s_{1},..., s_{n} \in \Omega_{\Gamma}$. 
Denote by ${\mathbb C}[z]_{d}$ the ${\mathbb C}$-vector space of 
polynomials over ${\mathbb C}$ of $z$ with degree $\leq d$. 
For $k > 1$, 
one has a ${\mathbb C}$-bilinear form $\Psi_{g,n;k}$ on 
$$
H^{0} \left( R_{\Gamma}, \Omega^{k}_{R_{\Gamma}} \left( k \sum_{j} s_{j} \right) \right) 
\times 
\left( H^{1} \left(\Gamma, {\mathbb C}[z]_{2k-2} \right) \oplus 
\left( {\mathbb C}[z]_{k-1} \right)^{\oplus n} \right) 
$$  
which is defined as 
$$
\Psi_{g,n;k} \left( \psi (dz)^{k}, \left( \xi, (f_{j})_{j} \right) \right) 
= \frac{1}{2 \pi \sqrt{-1}} \sum_{i=1}^{g} 
\oint_{\partial D_{i}} \psi \cdot \xi(\gamma_{i}) dz 
+ \sum_{j=1}^{n} {\rm Res}_{s_{j}} (\psi \cdot f_{j} dz). 
$$     
Then it is easy to see that $\Psi_{g,n;k}$ is non-degenerate since $\Psi_{g;k}$ is so. 
Let 
$$
\left\{ \xi_{1,n-1}, \xi_{2,1},..., \xi_{2,2k-2}, 
\xi_{i,0},..., \xi_{i,2k-2} \ (3 \leq l \leq g) \right\} 
$$
be the above basis of $H^{1} \left( \Gamma, {\mathbb C}[z]_{2k-2} \right)$. 
Then this basis together with 
$$
\left\{ \left. \left( z^{d_{1}},... , z^{d_{n}} \right) \ \right| \ 0 \leq d_{j} \leq k-1 \right\} 
$$
give a basis of 
$H^{1} \left( \Gamma, {\mathbb C}[z]_{2k-2} \right) \oplus 
\left( {\mathbb C}[z]_{k-1} \right)^{\oplus n}$. 
Therefore, there exists its dual basis of 
$H^{0} \left( R_{\Gamma}, \Omega^{k}_{R_{\Gamma}} \left( k \sum_{j} s_{j} \right) \right)$ 
for $\Psi_{g,n;k}$ which is also called {\it normalized}.  
\vspace{2ex}

5.3. {\it Explicit Riemann-Roch isomorphism}. 
\vspace{2ex}

{\sc Theorem 5.4}. 
\begin{it} 
The constant $c(g; k)$ in Proposition 5.3 becomes $\pm 1$. 
\end{it}
\vspace{2ex}

{\it Proof.} 
As in Section 4, 
let $C_{\Delta}$ be the generalized Tate curve over $A'_{\Delta}$ of  genus $g$ 
which is uniformized by $\Gamma_{\Delta}$, 
where the generators $\phi_{1},..., \phi_{g}$ of $\Gamma_{\Delta}$ 
is normalized as $x_{1} = 0$, $x_{-1} = \infty$, $x_{2} = 1$. 
As in stated in 4.1, 
one can also consider $C_{\Delta}$ as a family of Schottky uniformized 
Riemann surfaces by taking $x_{\pm i}$, $y_{i}$ 
as complex parameters $a_{\pm i}$, $q_{i}$ respectively 
such that $q_{i}$ are sufficiently small. 
By Proposition 4.2, 
$C_{\Delta}|_{y_{2} = \cdots = y_{g} = 0}$ becomes the stable curve 
$C'_{\Delta}$ obtained from the Tate curve 
${\cal E}_{y_{1}} = {\mathbb G}_{m}/\langle y_{1} \rangle$ over $A'_{g}[[y_{1}]]$ 
by identifying $x_{i} = x_{-i}$ $(2 \leq i \leq g)$. 
Hence by the properties of the relative dualizing sheaf (cf. [K, Section 1]), 
$H^{0} \left( C'_{\Delta}, \omega^{k}_{C'_{\Delta}/A'_{g}[[y_{1}]]} \right)$ 
becomes the subspace of 
$$
W = H^{0} \left( {\cal E}_{y_{1}}, \omega^{k}_{{\cal E}_{y_{1}}/A'_{g}[[y_{1}]]} 
\left( k \sum_{i=2}^{g} (x_{i} + x_{-i}) \right) \right)
$$ 
which consists of $\eta \in W$ satisfying 
${\rm Res}^{k}_{x_{i}} (\eta) = (-1)^{k} {\rm Res}^{k}_{x_{-i}} (\eta)$. 
Then by Proposition 4.1, there is a basis 
$\left\{ \eta_{1},..., \eta_{(2k-1)(g-1)} \right\}$ of the $A'_{g}[[y_{1}]]$-module
$H^{0} \left( C'_{\Delta}, \omega^{k}_{C'_{\Delta}/A'_{g}[[y_{1}]]} \right)$ 
which is normalized in the sense of 5.2. 
Since $k > 1$, 
$H^{1} \left( C, \omega_{C}^{k} \right) = \{ 0 \}$ 
for the dualizing sheaf $\omega_{C}$ on a stable curve $C$. 
Therefore, the natural homomorphism 
$$
H^{0} \left( C_{\Delta}, \omega^{k}_{C_{\Delta}/A'_{\Delta}} \right) \otimes_{A'_{\Delta}} 
\left( A'_{\Delta} / (y_{2},..., y_{g}) \right) \rightarrow 
H^{0} \left( C'_{\Delta}, \omega^{k}_{C'_{\Delta}/A'_{g}[[y_{1}]]} \right)
$$
is surjective, and hence there exists 
$$
\left\{ \theta_{1},..., \theta_{(2k-1)(g-1)} \right\} \subset 
H^{0} \left( C_{\Delta}, \omega^{k}_{C_{\Delta}/A'_{\Delta}} \right) 
$$ 
such that $\theta_{l}|_{y_{2} = \cdots y_{g} = 0} = \eta_{l}$ 
$\left( 1 \leq l \leq (2k-1)(g-1) \right)$. 

For each $i = 2,..., g$, 
let 
$$
z_{i} = \frac{(x_{i} - x_{-i})(z - x_{i})}{z - x_{-i}}, \   
z_{-i} = \frac{(x_{-i} - x_{i})(z - x_{-i})}{z - x_{i}} 
$$
be the local coordinates at $x_{i}$, $x_{-i}$ respectively 
such that 
$$
\lim_{z \rightarrow x_{i}} (z - x_{i})/z_{i} = 
\lim_{z \rightarrow x_{-i}} (z - x_{-i})/z_{-i} = 1. 
$$
Then the transformation matrices of 
$$
\left( z^{j} \right)_{0 \leq j \leq l} \mapsto 
\left( (z - x_{\pm i})^{j} \right)_{0 \leq j \leq l}, \ 
\left( (z - x_{\pm i})^{j} \right)_{0 \leq j \leq l} \mapsto 
\left( z_{\pm i}^{j} \ {\rm mod} \left( z_{\pm i}^{l+1} \right) 
\right)_{0 \leq j \leq l} 
$$
have determinant $1$, 
and hence we may replace the exterior product of 
the normalized basis with that of the basis of 
$H^{0} \left( C_{\Delta}, \omega^{k}_{C_{\Delta}/A'_{\Delta}} \right)$ 
dual to $\left( z_{\pm i}^{j} \right)_{1 \leq i \leq g}$. 
Since $C_{\Delta}$ is obtained as the deformation of ${\cal E}_{y_{1}}$ by the equation 
$z_{i} z_{-i} = -(x_{i} - x_{-i})^{2} y_{i}$ $(2 \leq i \leq g)$, 
$$
\frac{(dz_{i})^{k}}{z_{i}^{k+1}} \wedge \cdots \wedge 
\frac{(dz_{i})^{k}}{z_{i}^{2k-1}} 
= \pm \left( (x_{i} - x_{-i})^{2} y_{i} \right)^{-(k^{2}-k)/2} 
\frac{(dz_{-i})^{k}}{z_{-i}} \wedge \cdots \wedge 
\frac{(dz_{-i})^{k}}{z_{-i}^{k-1}}. 
$$
Let $\{ \xi_{1},..., \xi_{(2k-1)(g-1)} \}$ be the natural basis 
$$
\left\{ \xi_{1,k-1}, \xi_{2,1},..., \xi_{2,2k-2}, 
\xi_{i,0},..., \xi_{i,2k-2} \ (3 \leq i \leq g) \right\}
$$
given in 5.2. 
Then 
$$
\prod_{i=2}^{g} \left( (x_{i} - x_{-i})^{2} y_{i} \right)^{(k^{2}-k)/2} 
\det \left( \Psi_{g;k} \left( \theta_{l}, \xi_{m} \right) \right)_{l,m} 
$$ 
is a holomorphic function of $y_{2},..., y_{g}$ around the point defined as 
$y_{2} = \cdots = y_{g} = 0$. 
Furthermore, under $y_{i} \rightarrow 0$ $(2 \leq i \leq g)$, 
$\frac{1}{2 \pi \sqrt{-1}} \oint_{C_{\pm i}} \rightarrow {\rm Res}_{x_{\pm i}}$, 
and hence 
$$ 
\prod_{i=2}^{g} \left( (x_{i} - x_{-i})^{2} y_{i} \right)^{(k^{2}-k)/2} 
\det \left( \Psi_{g;k} \left( \theta_{l}, \xi_{m} \right) \right)_{l,m} 
\rightarrow \pm 1. 
$$ 

By definiton, the normalized basis $\left\{ \omega^{(k)}_{l} \right\}$ 
of holomorphic $k$-forms on $C_{\Delta}$ satisfies 
$$
\det \left( \Psi_{g;k} \left( \omega^{(k)}_{l}, \xi_{m} \right) \right)_{l,m} = 1 
$$ 
which implies that  
$$
\prod_{i=2}^{g} \left( (x_{i} - x_{-i})^{2} y_{i} \right)^{-(k^{2}-k)/2} 
\bigwedge_{l=1}^{(2k-1)(g-1)} \omega^{(k)}_{l}  
$$ 
gives a holomorphic section of $\omega_{C_{\Delta}}^{k}$ 
around $y_{2} = \cdots = y_{g} = 0$, 
and it becomes $\pm \bigwedge_{l=1}^{(2k-1)(g-1)} \eta_{l}$ 
at $y_{2} = \cdots = y_{g} = 0$. 
Then by applying Proposition 3.3, Theorem 3.2 to Proposition 5.3, 
and by using Propositions 4.2, 4.3, we have 
$$
c(g;k) \cdot \mu_{1,2g-2;k} \left( \bigwedge_{l=1}^{(2k-1)(g-1)} \eta_{l} \right) = 
\pm \prod_{m=1}^{\infty} \left( 1 - y_{1}^{m} \right)^{2(d_{k}-1)} 
\cdot \prod_{i=2}^{g} (x_{i} - x_{-i})^{-(k^{2}-k)} 
$$  
under the trivialization by the basis of 
$H^{0} \left( C'_{\Delta}, \omega^{k}_{C'_{\Delta}/A'_{g}[[y_{1}]]} \right)$ and 
$d(z - x_{\pm i})$ $(2 \leq i \leq g)$. 
Since $\mu_{1,2g-2;k} \left( \bigwedge_{l=1}^{(2k-1)(g-1)} \eta_{l} \right)$ and 
the right hand side is primitive, 
the constant $c(g;k)$ is seen to be $\pm 1$. 
\ $\square$ 
\vspace{2ex}

{\sc Theorem 5.5}. 
\begin{it} 
The Riemann-Roch isomorphism 
$$
\lambda_{g;k}^{\otimes 12} \otimes \delta_{g}^{\otimes (-1)} 
\stackrel{\sim}{\rightarrow}  
\kappa_{g}^{\otimes d_{k}}
$$
given by Theorem 3.1 maps $\left( F_{k} \omega^{(k)} \right)^{12}$ to 
$\pm \left( {\rm ACS}_{g} \right)^{-d_{k}}$. 
\end{it} 
\vspace{2ex}

{\it Proof.} 
First, we prove the assertion when $k = 1$. 
By the definition of ${\rm ACS}_{g}$, 
the image of $(F_{1} \omega)^{12}$ is a constant multiple of 
$\left( {\rm ACS}_{g} \right)^{-1}$, 
and by the arithmetic Schottky uniformization theory, 
the both are primitive. 
Therefore, they are equal up to a sign. 
Second, we assume that $k > 1$. 
Since the norm of ${\rm ACS}_{g}$ is a constant multiple of 
$\exp \left( -S_{\rm L}/(2 \pi) \right)$, 
Theorem 5.1 implies that the image of $\left( F_{k} \omega^{(k)} \right)^{12}$ 
is a certain multiple of $({\rm ACS}_{g})^{-d_{k}}$. 
Since $F_{1}$ and $F_{k}$ are expressed as primitive elements of $A'_{\Delta}$, 
by Theorem 5.4, $\omega^{(k)}$ gives a trivialization of $\lambda_{g;k}$ and is primitive. 
Therefore, the image of $\left( F_{k} \omega^{(k)} \right)^{12}$ is also primitive. 
This completes the proof. 
\ $\square$ 
\vspace{2ex} 

{\sc Theorem 5.6}. 
\begin{it} 
Let $\left\{ \psi_{l} \right\}$ be the normalized basis of 
$H^{0} \left( R_{\Gamma}, \Omega^{k}_{R_{\Gamma}} \left( k \sum_{j} s_{j} \right) \right)$ 
for Schottky uniformized Riemann surfaces $R_{\Gamma}$ with marked points 
$s_{1},..., s_{n}$ 
Then 
\end{it} 
$$
\mu_{g,n;k} 
\left( \bigwedge_{l} \psi_{l} \right) 
= \pm \frac{F_{1}^{d_{k}}}{F_{k}} \omega^{\otimes d_{k}} 
\otimes \bigotimes_{j=1}^{n} d(z - s_{j})^{\otimes (k^{2}-k)/2}. 
$$

{\it Proof.} 
Let $\left( \pi : C \rightarrow S; \sigma_{1},..., \sigma_{n} \right)$ 
be an $n$-pointed stable curve of genus $g > 1$, 
Then for non-negative integers $l_{1},..., l_{n} \leq k$, 
the $l_{i}$th residue map gives the exact sequence 
\begin{eqnarray*} 
0 
& \rightarrow & 
\pi_{*} \left( \omega_{C/S}^{k} \left( \sum_{j=1}^{n} l_{j} \sigma_{j} \right) \right) 
\rightarrow  
\pi_{*} \left( \omega_{C/S}^{k} \left( \sum_{j \neq i} l_{j} \sigma_{j} + 
(l_{i} + 1) \sigma_{i} \right) \right) 
\\ 
& \rightarrow & 
{\cal O}_{S} (d z_{\sigma_{i}})^{k-l_{i}-1} \rightarrow 0, 
\end{eqnarray*} 
where $d z_{\sigma_{i}}$ denotes the local coordinate around $\sigma_{i}$. 
Therefore, by induction, 
$$
\det \pi_{*} \left( \omega_{C/S}^{k} \left( k \sum_{j=1}^{n} \sigma_{j} 
\right) \right) 
\cong 
\det \pi_{*} \left( \omega_{C/S}^{k} \right) \cdot 
\prod_{j=1}^{n} (d z_{\sigma_{i}})^{(k^{2} - k)/2}. 
$$
If $C$ consists of Schottky uniformized Riemann surfaces 
with marked points $s_{1},..., s_{n}$, 
then the residue map sends the normalized basis to 
$\left\{ d(z - s_{i})^{k - l_{i} - 1} \right\}$, 
and hence the assertion follows from Theorem 5.4. 
\ $\square$ 
\vspace{2ex}

5.4. {\it Completed formulae of Zograf and McIntyre-Takhtajan}. 
\vspace{2ex}

{\sc Theorem 5.7}. 
\begin{it}
Let $c_{g}$ be the constant given in Theorem 2.1. 
Then 
$$
c_{g} = (2 \pi)^{2g} \cdot \exp \left( 
\frac{(g-1) \left( 24 \zeta'_{\mathbb Q}(-1) + 1 \right)}{6} \right). 
$$ 
\end{it}

{\it Proof.} 
Let ${\cal C}$ be an algebraic curve of genus $g$ over a scheme $S$, 
and represent $\Omega_{{\cal C}/S}$ as ${\cal O}_{\cal C}({\cal D})$ and 
${\cal O}_{\cal C}({\cal D}')$ such that the supports of ${\cal D}$, ${\cal D}'$ 
are disjoint. 
Then by definition, $\langle \mbox{\boldmath $1$}, \mbox{\boldmath $1$} \rangle$ 
is a section of the line bundle 
$\left\langle \Omega_{{\cal C}/S}, \Omega_{{\cal C}/S} \right\rangle$ on $S$ 
given by the Deligne pairing.  
Denote by $S_{\rm A}[\log \rho]$ the Liouville action 
defined in [A, Definition 4.1]. 
Then by the Gauss-Bonnet theorem, 
$S_{\rm L} = 2 \pi S_{\rm A} [\log \rho]$ 
when $\rho$ is the hyperbolic metric with constant curvature $-1$. 
Since the tangent line bundle ${\cal T}_{{\cal C}/S}$ of ${\cal C}/S$ is 
$\Omega_{{\cal C}/S}^{\otimes (-1)}$, 
$\left\langle {\cal T}_{{\cal C}/S}, {\cal T}_{{\cal C}/S} \right\rangle = 
\left\langle \Omega_{{\cal C}/S}, \Omega_{{\cal C}/S} \right\rangle$. 
Furthermore, if ${\cal C}/S$ is a family of Riemann surfaces, 
then [A, Corollary 5.2] implies that  
$$
\exp \left( S_{\rm L}/({2 \pi}) \right) = 
\left\| \langle \mbox{\boldmath $1$}, \mbox{\boldmath $1$} \rangle \right\| 
\cdot e^{2g-2}. 
$$
Let $\omega_{i}$ $(1 \leq i \leq g)$ be as in 4.2 which becomes the normalized basis of 
holomorphic $1$-forms on Schottky uniformized Riemann surfaces $R_{\Gamma}$, 
and put $\omega = \omega_{1} \wedge \cdots \wedge \omega_{g}$. 
Then $\int_{\alpha_{i}} \omega_{j} = 2 \pi \sqrt{-1} \delta_{ij}$, 
and hence 
$$
\frac{\det {\rm Im}(\tau)}{\det' \Delta_{h}} = 
\frac{\left\| \omega/(2 \pi \sqrt{-1})^{g} \right\|^{2}_{\rm H}}{\det' \Delta_{h}} = 
\frac{\left\| \omega \right\|^{2}_{\rm Q}}{(2 \pi)^{2g}}. 
$$
Furthermore, by Theorem 2.1, there exists an isometry 
$\lambda_{1}^{\otimes 12} \stackrel{\sim}{\rightarrow} 
\left\langle \Omega_{{\cal C}/S}, \Omega_{{\cal C}/S} \right\rangle$ 
which sends $(F_{1} \omega)^{12}$ to 
$\langle \mbox{\boldmath $1$}, \mbox{\boldmath $1$} \rangle$ 
up to a nonzero constant. 
Since $\omega$ and $\langle \mbox{\boldmath $1$}, \mbox{\boldmath $1$} \rangle$ 
give primitive local sections of $\lambda_{g;1}$ and $\kappa_{g}$ respectively, 
the Riemann-Roch isomorphism sends $(F_{1} \omega)^{12}$ to 
$\pm \langle \mbox{\boldmath $1$}, \mbox{\boldmath $1$} \rangle$. 
Therefore, by Theorem 3.1, 
\begin{eqnarray*}
\left( \frac{\det {\rm Im}(\tau)}{\det' \Delta_{h}} \right)^{6} 
& = & 
\left\| \omega \right\|_{\rm Q}^{12} \cdot (2 \pi)^{-12g} 
\\ 
& = & 
\left\| \left\langle \mbox{\boldmath $1$}, \mbox{\boldmath $1$} \right\rangle \right\|  \cdot |F_{1}|^{-12} \cdot e^{a(g)} \cdot (2 \pi)^{-12g} 
\\ 
& = & 
\exp \left( \frac{S_{\rm L}}{2 \pi} \right) \cdot |F_{1}|^{-12} \cdot 
\exp \left( (1-g) \left( 24 \zeta'_{\mathbb Q}(-1) + 1 \right) \right) 
\cdot (2 \pi)^{-12g}. 
\end{eqnarray*}
Comparing this equality with the definition of $c_{g}$ given in Theorem 2.1, 
we have  
$$
c_{g} = (2 \pi)^{2g} \cdot \exp \left( 
\frac{(g-1) \left( 24 \zeta'_{\mathbb Q}(-1) + 1 \right)}{6} \right) 
$$ 
which completes the proof. 
\ $\square$ 
\vspace{2ex}

{\sc Theorem 5.8}. 
\begin{it}
The constant $c_{g;k}$ in Theorem 5.2 is determined as  
$$
c_{g;k} = \exp \left( \frac{(g-1) \left( 24 \zeta'_{\mathbb Q}(-1) + 2 d_{k} - 1 \right)}{6} \right). 
$$ 
\end{it}

{\it Proof.}  
By Theorem 3.1, 
\begin{eqnarray*}
\left\| \omega^{(k)} \right\|_{\rm Q}^{12}  
& = & 
\left\| \langle \mbox{\boldmath $1$}, \mbox{\boldmath $1$} \rangle \right\|^{d_{k}} \cdot 
|F_{k}|^{-12} \cdot e^{a(g)} 
\\ 
& = & 
\exp \left( \frac{S_{\rm L}}{2 \pi} \right)^{d_{k}} \cdot |F_{k}|^{-12} \cdot 
\exp \left( (1-g) \left( 24 \zeta'_{\mathbb Q}(-1) + 2 d_{k} - 1 \right) \right) 
\end{eqnarray*}
which completes the proof. 
\ $\square$ 
\vspace{2ex}

\begin{center}
6. {\bf Rationality of Ruelle zeta values} 
\end{center}

6.1. {\it Mumford isomorphism and discriminant}. 
Assume that $g > 1$, 
let $\left( \Gamma; \gamma_{1},..., \gamma_{g} \right)$ be 
a marked normalized Schottky group, 
and $F_{k}(\Gamma; \gamma_{1},..., \gamma_{g})$ 
be the value of $F_{k}$ at the corresponding point on ${\mathfrak S}_{g}$. 
Recall that $R_{\Gamma}$ denotes the Riemann surface uniformized by $\Gamma$ 
which is equivalently the algebraic curve over ${\mathbb C}$ associated with $\Gamma$. 
Assume that there exists a sub ${\mathbb Z}$-algebra $R$ of ${\mathbb C}$ 
over which $R_{\Gamma}$ has a model $C_{\Gamma}$ as a stable curve, 
and that there exist $R$-basis $\{ u_{1},..., u_{g} \}$ of 
$H^{0} \left( C_{\Gamma}, \omega_{C_{\Gamma}/R} \right)$ 
and $\{ v_{1},..., v_{(2k-1)(g-1)} \}$ of 
$H^{0} \left( C_{\Gamma}, \omega^{k}_{C_{\Gamma}/R} \right)$ for $k > 1$. 
Then their periods are defined as 
$$
\Omega_{1} = \det \left( \frac{1}{2 \pi \sqrt{-1}} 
\oint_{C_{i}} u_{j} \right)_{1 \leq i, j \leq g}, 
$$ 
and  
$$
\Omega_{k} = \det \left( \frac{1}{2 \pi \sqrt{-1}} \sum_{i=1}^{g} 
\oint_{C_{i}} v_{l} \cdot \xi_{m}(\gamma_{i}) 
\right)_{1 \leq l,m \leq (2k-1)(g-1)}, 
$$
where $C_{1},..., C_{g} \subset {\mathbb P}^{1}({\mathbb C})$ are given in 2.2, 
and $\left\{ \xi_{1},..., \xi_{(2k-1)(g-1)} \right\}$ is the basis 
$$
\left\{ \xi_{1,k-1}, \ \xi_{2,1},..., \xi_{2,2k-2}, \ 
\xi_{i,0},..., \xi_{i,2k-2} \ (3 \leq i \leq g) \right\}
$$ 
of $H^{1} \left( \Gamma, {\mathbb C}[z]_{2k-2} \right)$ given in 5.2. 

Following [S], let $D(C_{\Gamma})$ denote the discriminant ideal 
associated with $C_{\Gamma}$, 
which is defined as the ideal of $R$ corresponding to the boundary 
$\partial {\cal M}_{g} = \overline{\cal M}_{g} - {\cal M}_{g}$ of $\overline{\cal M}_{g}$ 
via the morphism ${\rm Spec} (R) \rightarrow \overline{\cal M}_{g}$ 
associated with $C_{\Gamma}$. 
Then we consider the rationality of  
$$
\frac{F_{k}(\Gamma; \gamma_{1},..., \gamma_{g})}{\Omega_{k}}, 
$$
and its relation with $D(C_{\Gamma})$. 
\vspace{2ex}

{\sc Theorem 6.1}. 
\begin{it} 
\begin{itemize}

\item[{\rm (1)}] 
Under the above notations,  
the ratio 
$$
\frac{\left( F_{1}(\Gamma; \gamma_{1},..., \gamma_{g}) /\Omega_{1} \right)^{d_{k}}}
{F_{k}(\Gamma; \gamma_{1},..., \gamma_{g}) /\Omega_{k}} 
$$
belongs to the power $D(C_{\Gamma})^{(k^{2}-k)/2}$ of $D(C_{\Gamma})$. 

\item[{\rm (2)}] 
If $R$ is a discrete valuation ring and $d(C_{\Gamma})$ is a generator of $D(C_{\Gamma})$, 
then 
$$
\frac{\left( F_{1}(\Gamma; \gamma_{1},..., \gamma_{g}) /\Omega_{1} \right)^{d_{k}}}
{F_{k}(\Gamma; \gamma_{1},..., \gamma_{g}) /\Omega_{k}} 
\in d(C_{\Gamma})^{(k^{2}-k)/2} \cdot R^{\times}, 
$$
where $R^{\times}$ denotes the unit group of $R$. 

\item[{\rm (3)}] 
If $C_{\Gamma}$ is smooth over $R$, 
then 
$$
\frac{\left( F_{1}(\Gamma; \gamma_{1},..., \gamma_{g}) /\Omega_{1} \right)^{d_{k}}}
{F_{k}(\Gamma; \gamma_{1},..., \gamma_{g}) /\Omega_{k}} 
\in R^{\times}.
$$
\end{itemize}
\end{it}

{\it Proof.} 
By Theorem 3.2, 
$$
\frac{\left( F_{1}(\Gamma; \gamma_{1},..., \gamma_{g}) /\Omega_{1} \right)^{d_{k}}}
{F_{k}(\Gamma; \gamma_{1},..., \gamma_{g}) /\Omega_{k}} 
$$
gives (up to a sign) the evaluation of the Mumford isomorphism $\mu_{g,0;k}$ 
on $C_{\Gamma}$ under the trivializations of $\lambda_{g;1}$ and $\lambda_{g;k}$ 
by $u_{1} \wedge \cdots \wedge u_{g}$ and 
$v_{1} \wedge \cdots \wedge v_{(2k-1)(g-1)}$ respectively. 
Hence the assertions follows from Theorem 5.4. 
\ $\square$ 
\vspace{2ex}

6.2. {\it Rationality of Ruelle zeta values}. 
Let $\Gamma$ be a Schottky group, 
and $X_{\Gamma} = {\mathbb H}^{3}/\Gamma$ be the hyperbolic $3$-manifold 
uniformized by $\Gamma$. 
Then for each $\gamma \in \Gamma$, $-\log |q_{\gamma}|$ is the length 
of the closed geodesic on $X_{\Gamma}$ corresponding to $\gamma$, 
and hence the Ruelle zeta function of $X_{\Gamma}$ becomes 
$$
Z_{\Gamma}(s) = \prod_{\{ \gamma \}} \left( 1 - |q_{\gamma}|^{s} \right)^{-1}, 
$$
where $\{ \gamma \}$ runs over primitive conjugacy classes of $\Gamma$. 
It is known (cf. [MT, 5.2]) that $Z_{\Gamma}(s)$ is absolutely convergent 
if ${\rm Re}(s) \geq 2$. 
We assume that a marked normalized Schottky group 
$( \Gamma; \gamma_{1},..., \gamma_{g})$ is contained in 
$PSL_{2}({\mathbb R})$, 
and apply Theorem 6.1 to showing the rationality of 
the {\it modified} Ruelle zeta values 
$$
\widetilde{Z}_{\Gamma}(k) = Z_{\Gamma}(k) 
\frac{\left( 1 - q_{\gamma_{1}}^{k} \right)^{2} 
\left( 1 - q_{\gamma_{2}}^{k} \right)}
{\left( 1 - q_{\gamma_{2}}^{k-1} \right)} 
$$
for integers $k > 1$. 
\vspace{2ex}

{\sc Theorem 6.2}. 
\begin{it} 
Assume that a Schottky group $\Gamma$ is contained in $PSL_{2}({\mathbb R})$, 
and that $R_{\Gamma}$ has a model as a stable curve $C_{\Gamma}$ 
over a sub ${\mathbb Z}$-algebra $R$ of ${\mathbb C}$. 
Let $k$ is an integer $> 1$, and put 
$$
c(\Gamma) = \frac{F_{1}(\Gamma; \gamma_{1},..., \gamma_{g})}{\Omega_{1}}. 
$$ 
\begin{itemize}

\item[{\rm (1)}] 
If $K$ denotes the quotient field of $R$, 
then 
$$
\widetilde{Z}_{\Gamma}(k) 
\in 
\frac{\Omega_{k+1}}{\Omega_{k}} \cdot c(\Gamma)^{12k} \cdot K^{\times}. 
$$

\item[{\rm (2)}] 
If $R$ is a discrete valuation ring, 
then 
$$
\widetilde{Z}_{\Gamma}(k) 
\in 
\frac{\Omega_{k+1}}{\Omega_{k}} \cdot c(\Gamma)^{12k} \cdot d(C_{\Gamma})^{-k} 
\cdot R^{\times}. 
$$

\item[{\rm (3)}] 
If $C_{\Gamma}$ is smooth over $R$, 
then 
$$
\widetilde{Z}_{\Gamma}(k) 
\in 
\frac{\Omega_{k+1}}{\Omega_{k}} \cdot c(\Gamma)^{12k} \cdot R^{\times}. 
$$
\end{itemize}
\end{it}

{\it Proof.} 
By assumption, 
$q_{\gamma} = |q_{\gamma}|$ for any $\gamma \in \Gamma - \{ 1 \}$, 
and hence 
$$
\widetilde{Z}_{\Gamma}(k) = 
\frac{F_{k+1}(\Gamma; \gamma_{1},..., \gamma_{g})}
{F_{k}(\Gamma; \gamma_{1},..., \gamma_{g})}. 
$$
Therefore, the assertions follows from Theorem 6.1. 
\ $\square$ 
\vspace{2ex}

{\sc Corollary 6.3}.  
\begin{it} 
Let the notation and assumption be as in Theorem 6.2, 
and assume that $R$ is a Dedekind ring with quotient field $K$. 
For $j \in \{ 1, k, k+1 \}$, 
let $\left\{ v_{i}^{(j)} \right\}_{i}$ be a $K$-basis of 
$H^{0} \left( C_{\Gamma}, \omega^{j}_{C_{\Gamma}/R} \right) \otimes_{R} K$, 
and denote by $\Omega_{j}$ the period of this basis. 
\begin{itemize}

\item[{\rm (1)}] 
There exists an element $r(C_{\Gamma})$ of $K^{\times}$ such that 
$$
\widetilde{Z}_{\Gamma}(k) 
\in 
\frac{\Omega_{k+1}}{\Omega_{k}} \cdot c(\Gamma)^{12k} \cdot r(C_{\Gamma}) 
\cdot R^{\times}. 
$$

\item[{\rm (2)}] 
For each prime ideal $P$ of $R$, 
${\rm ord}_{P} \left( r(C_{\Gamma}) \right)$ is given by 
$$
k \cdot {\rm ord}_{P} \left( d(C_{\Gamma}) \right) 
+ 12 k \cdot {\rm ord}_{P} \left( \bigwedge v_{i}^{(1)} \right) 
+ {\rm ord}_{P} \left( \bigwedge v_{i}^{(k)} \right) 
- {\rm ord}_{P} \left( \bigwedge v_{i}^{(k+1)} \right), 
$$
where ${\rm ord}_{P} \left( \bigwedge v_{i}^{(j)} \right)$ denotes 
the order at $P$ of the fractional ideal of $R_{P}$ generated by 
the exterior product of $\left\{ v_{i}^{(j)} \right\}_{i}$ under an identification 
$$
R_{P} = 
\det \left( H^{0} \left( C_{\Gamma}, \omega^{j}_{C_{\Gamma}/R} \right) \right) 
\otimes_{R} R_{P} 
$$
as $R_{P}$-modules. 
\end{itemize}
\end{it}

\begin{center} 
{\bf References}
\end{center}

\begin{itemize}

\item[{[A]}] 
E. Aldrovandi, 
On hermitian-holomorphic classes related to uniformization, 
the dilogarithm, and the Liouville action, 
{\it Comm. Math. Phys.} {\bf 251} (2004), 27--64. 

\item[{[B]}] 
A. Beilinson, 
Higher regulators and values of $L$-functions, 
{\it J. Sov. Math.} {\bf 30} (1985), 2036--2070. 

\item[{[BK]}] 
S. Bloch and K. Kato, 
$L$-functions and Tamagawa numbers of motives, 
in {\it The Grothendieck Festschrift Vol. I, Progr. Math.} Vol. 86, 
Birkh\"{a}user, Boston, 1990, pp. 333--400. 

\item[{[CS]}] 
S.-S. Chern and J. Simons, 
Characteristic forms and geometric invariants, 
{\it Ann. Math.} {\bf 99} (1974), 48--69. 

\item[{[D1]}] 
P. Deligne, 
Valeurs de fonctions $L$ et periodes d'integrales, 
in {\it Automorphic forms, representations, and $L$-functions, 
Proc. Symp. in Pure Math.} Vol. 33, Part 2, 
Amer. Math. Soc., 1979, pp. 313--346. 

\item[{[D2]}] 
P. Deligne, 
Le det\'{e}rminant de la cohomologie, 
in {\it Current trends in arithmetical algebraic geometry, Contemp. Math.} Vol. 67, 
Amer. Math. Soc., 1987, pp. 93--177. 

\item[{[DM]}] 
P. Deligne and D. Mumford, 
The irreducibility of the space of curves of given genus, 
{\it Publ. Math. IHES} \ {\bf 36} (1969), 75--109. 

\item[{[F]}] 
D. S. Freed, 
Classical Chern-Simons theory, Part 1, 
{\it Adv. Math.} {\bf 113} (1995), 237--303. 

\item[{[Fr]}] 
G. Freixas i Montplet, 
An arithmetic Riemann-Roch theorem for pointed stable curves, 
{\it Ann. Scient. \'{E}c. Norm. Sup.} {\bf 42} (2009), 335--369. 

\item[{[Fri]}] 
D. Fried, 
Analytic torsion and closed geodesics on hyperbolic manifolds, 
{\it Invent. Math.} {\bf 84} (1986), 523--540. 

\item[{[GP]}] 
Y. Gon and J. Park, 
The zeta functions of Ruelle and Selberg for hyperbolic manifolds with cusps, 
{\it Math. Ann.} {\bf 346} (2010), 719--767. 

\item[{[GS]}] 
H. Gillet and Soul\'{e}, 
Arithmetic intersection theory, 
{\it Publ. Math. IHES} \ {\bf 72} (1990), 94--174. 

\item[{[GM]}] 
C. Guillarmou and S. Moroianu, 
Chern-Simons line bundle on Teichm\"{u}ller space, 
{\it Geometry \& Topology} {\bf 18} (2014), 327--377. 

\item[{[I1]}] 
T. Ichikawa, 
$P$-adic theta functions and solutions of the KP hierarchy, 
{\it Comm. Math. Phys.} {\bf 176} (1996), 383--399. 

\item[{[I2]}] 
T. Ichikawa, 
Generalized Tate curve and integral Teichm\"{u}ller modular forms, 
{\it Amer. J. Math.} {\bf 122} (2000), 1139--1174. 

\item[{[I3]}] 
T. Ichikawa, 
Universal periods of hyperelliptic curves and their applications, 
{\it J. Pure Appl. Algebra} {\bf 163} (2001), 277--288. 

\item[{[K]}] 
F. F. Knudsen, 
The projectivity of the moduli space of stable curves II, III, 
{\it Math. Scand.} {\bf 52} (1983), 161--199, 200--212. 

\item[{[KK1]}] 
A. Kokotov and D. Korotkin, 
Tau-functions on Hurwitz spaces, 
{\it Math. Phys. Anal. Geom.} {\bf 7} (2004), 47--96. 

\item[{[KK2]}] 
A. Kokotov and D. Korotkin, 
Tau-functions on spaces of abelian differentials 
and higher genus generalizations of Ray-Singer formula, 
{\it J. Diff. Geom.} {\bf 82} (2009), 35--100. 

\item[{[Kr]}] 
K. Krasnov, 
Holography and Riemann surfaces, 
{\it Adv. Theor. Math. Phys.} {\bf 4} (2000), 929--979. 

\item[{[MD]}] 
Y. Manin and V. Drinfeld, 
Periods of $p$-adic Schottky groups, 
{\it J. Reine Angew. Math.} {\bf 262/263} (1973), 239--247. 

\item[{[MP]}] 
A. McIntyre and J. Park, 
Tau functions and Chern-Simons invariant, 
{\it Adv. Math.} {\bf 262} (2014), 1--58. 

\item[{[MT]}] 
A. McIntyre and L. Takhtajan, 
Holomorphic factorization of determinants of laplacians on Riemann surfaces 
and a higher genus generalization of Kronecker's first limit formula, 
{\it GAFA} {\bf 16} (2006), 1291--1323. 

\item[{[MoT]}] 
M. Morishita and Y. Terashima, 
Chern-Simons variation and Deligne cohomology, 
in {\it Spectral Analysis in Geometry and Number Theory, 
Contemporary Math.} {\bf 484} (2009), 127--134. 

\item[{[Mu1]}] 
D. Mumford, 
An analytic construction of degenerating curves over complete local rings, 
{\it Compositio Math.} {\bf 24} (1972), 129--174. 

\item[{[Mu2]}] 
D. Mumford, 
Stability of projective varieties, 
{\it L'Ens. Math.} {\bf 23} (1977), 39--110. 

\item[{[P]}] 
J. Park, 
Analytic torsion and Ruelle zeta functions for hyperbolic manifolds with cusps, 
{\it J. Funct. Anal.} {\bf 257} (2009), 1713--1758. 

\item[{[Q]}] 
D. Quillen, 
Determinants of Cauchy-Riemann operators over a Riemann surface, 
{\it Funct. Anal. Appl.} {\bf 19} (1986), 31--34. 

\item[{[RSW]}] 
T. R. Ramadas, I. M. Singer and J. Weitsman, 
Some comments on Chern-Simons gauge theory, 
{\it Commun. Math. Phys.} {\bf 126} (1989), 409--420. 

\item[{[RS]}] 
D. B. Ray and I. M. Singer, 
$R$-Torsion and the Laplacian on Riemannian manifolds, 
{\it Adv. Math.} {\bf 7} (1971), 145--210. 

\item[{[S]}] 
T. Saito, 
Conductor, discriminant, and the Noether formula of arithmetic surfaces, 
{\it Duke Math. J.} {\bf 57} (1988), 151--173. 

\item[{[Si]}] 
J. H. Silverman, 
{\it Advanced topics in the arithmetic of elliptic curves, 
Graduate Texts in Math.} Vol. 151, Springer, 1994. 

\item[{[Su1]}] 
K. Sugiyama, 
An analog of the Iwasawa conjecture for a compact hyperbolic threefold, 
{\it J. Reine Angew. Math.} {\bf 613} (2007), 35--50. 

\item[{[Su2]}] 
K. Sugiyama, 
A special value of Ruelle $L$-function and the theorem of 
Cheeger and M\"{u}ller, 
arXiv:0803.2079v1 

\item[{[Su3]}] 
K. Sugiyama, 
The Taylor expansion of Ruelle $L$-function at the origin 
and the Borel regulator, 
arXiv:0804.2715v1 

\item[{[Su4]}] 
K. Sugiyama, 
On geometric analogues of the Birch and Swinnerton-Dyer conjecture 
for low dimensional hyperbolic manifolds, 
in {\it Spectral analysis in geometry and number theory, Contemp. Math.} 
Vol. 484, Amer. Math. Soc., 2009. pp. 267--286. 

\item[{[Su5]}] 
K. Sugiyama, 
On geometric analogues of Iwasawa main conjecture for a hyperbolic threefold, 
in {\it Noncommutativity and singularities, Adv. Stud. Pure Math.} 
Vol. 55, Math. Soc. Japan, Tokyo, 2009, pp. 117--135. 

\item[{[TT]}] 
L. A. Takhtajan and L.-P. Teo, 
Liouville action and Weil-Petersson metric on deformation spaces, 
global Kleinian reciprocity and holography, 
{\it Comm. Math. Phys.} {\bf 239} (2003), 183--240. 

\item[{[T]}] 
J. Tate, 
A review of non-archimedean elliptic functions, 
in {\it Elliptic Curves, Modular forms, \& Fermat's Last Theorem,} 
International Press, Boston, 1995, pp. 162--184. 

\item[{[W]}] 
L. Weng, 
$\Omega$-admissible theory II, 
Deligne pairings over moduli spaces of punctured Riemann surfaces, 
{\it Math. Ann.} {\bf 320} (2001), 239--283. 

\item[{[Y]}] 
T. Yoshida, 
The $\eta$-invariant of hyperbolic $3$-manifolds, 
{\it Invent. Math.} {\bf 81} (1985), 473--514.

\item[{[Z1]}] 
P. G. Zograf, 
Liouville action on moduli spaces and uniformization 
of degenerate Riemann surfaces, 
{\it Algebra, i Analiz} {\bf 1} (1989), 136--160 (Russian), 
English translation in {\it Leningrad Math. J.} {\bf 1} (1990), 941--965. 

\item[{[Z2]}] 
P. G. Zograf, 
Determinants of Laplacians, Liouville action, 
and an analogue of the Dedekind $\eta$-function on Teichm\"{u}uller space, 
Unpublished manuscript (1997). 

\item[{[ZT]}] 
P. G. Zograf and L. A. Takhtajan, 
On the uniformization of Riemann surfaces and on the Weil-Petersson metric 
on the Teichm\"{u}ller and Schottky spaces, 
{\it Math. USSR-Sb.} {\bf 60} (1988), 297--313. 

\end{itemize}

\begin{flushleft}
{\sc Department of Mathematics, 
Graduate School of Science and Engineering, 
Saga University, Saga 840-8502, Japan} 
\\
{\it E-mail}: ichikawa@ms.saga-u.ac.jp 
\end{flushleft} 

\end{document}